\documentclass[a4paper,12pt]{article}

\usepackage{amsmath,amssymb,latexsym,amsthm,enumerate,nicefrac,bbm,hyperref,xcolor,graphicx}
\usepackage{cite}
\usepackage[left=20mm,right=18mm,top=1.5cm,bottom=1cm,includeheadfoot]{geometry}

\usepackage[T1]{fontenc}
\usepackage[utf8]{inputenc}
\usepackage[english]{babel}
\usepackage[capitalise,compress]{cleveref}

\usepackage{csquotes}

\newtheorem{theo}{Theorem}[section]

\newtheorem{lemma}[theo]{Lemma}
\newtheorem{propo}[theo]{Proposition}
\newtheorem{corollary}[theo]{Corollary}
\newtheorem{setting}[theo]{Setting}

\theoremstyle{definition}

\newtheorem{defi}[theo]{Definition}

\newcommand*{\wt}{\widetilde}

\newcommand*{\N}{\mathbb{N}}
\newcommand*{\R}{\mathbb{R}}
\newcommand*{\Z}{\mathbb{Z}}

\newcommand*{\E}{\mathbb{E}}

\newcommand*{\X}{\mathbb{X}}
\newcommand*{\Y}{\mathbb{Y}}

\newcommand*{\Dm}{\mathfrak{D}}
\newcommand*{\dm}{d}

\newcommand*{\w}{\mathfrak{w}}  %weight function
\newcommand*{\W}{\mathfrak{W}}  %"weighted" $L^\infty$ 
\newcommand*{\A}{\widetilde{\mathcal{A}}}   %sigma alg on R^A

\newcommand*{\var}{\operatorname{Var}}

\newcommand{\dkl}[3]{\big(#1(#2)\big)(#3)}
\newcommand{\pset}[1]{\mathnormal{2}^{#1}}

\newenvironment{mproof}[1]{\noindent\textit{Proof of {#1}.}}{\hfill \qed}
\title{Nonlinear Monte Carlo methods with polynomial\\ runtime for Bellman equations of discrete time high-dimensional stochastic optimal control problems}

\author{Christian Beck$^{1}$, Arnulf Jentzen$^{2,3}$, \\
Konrad Kleinberg$^{4}$, and Thomas Kruse$^{5}$\bigskip\\
\small{$^1$ Department of Mathematics, ETH Zurich,}\\
\small{Switzerland; e-mail: 
\texttt{christian.beck}\textcircled{\texttt{a}}\texttt{yahoo.de}}\\
\small{$^2$ Applied Mathematics: Institute for Analysis and Numerics, }\\ 
\small{University of M\"unster, Germany; e-mail: \texttt{ajentzen}\textcircled{\texttt{a}}\texttt{uni-muenster.de}}\\
\small{$^3$ School of Data Science and Shenzhen Research Institute of Big Data,}\\
\small{The Chinese University of Hong Kong, Shenzhen,}\\
\small{China; e-mail: \texttt{ajentzen}\textcircled{\texttt{a}}\texttt{cuhk.edu.ch}}\\
\small{$^4$ Department of Mathematics \& Informatics, University of Wuppertal,}\\
\small{Germany; e-mail: \texttt{kleinberg}\textcircled{\texttt{a}}\texttt{uni-wuppertal.de}}\\
\small{$^5$ Department of Mathematics \& Informatics, University of Wuppertal,}\\
\small{Germany; e-mail: \texttt{tkruse}\textcircled{\texttt{a}}\texttt{uni-wuppertal.de}}
}

\begin{document}
\maketitle

\begin{abstract}
Discrete time \emph{stochastic optimal control} problems and \emph{Markov decision processes} (MDPs), respectively, serve as fundamental models for problems that involve sequential decision making under uncertainty and as such constitute the theoretical foundation of \emph{reinforcement learning}. In this article we study the numerical approximation of MDPs with infinite time horizon, finite control set, and general state spaces. Our set-up in particular covers infinite-horizon optimal stopping problems of discrete time Markov processes. A key tool to solve MDPs are \emph{Bellman equations} which characterize the value functions of the MDPs and determine the optimal control strategies.
By combining ideas from the \emph{full-history recursive multilevel Picard approximation method}, which was recently introduced to solve certain nonlinear partial differential equations, and ideas from \emph{$Q$-learning} we introduce a class of suitable \emph{nonlinear Monte Carlo methods} and prove that the proposed methods do overcome the \emph{curse of dimensionality} in the numerical approximation of the solutions of Bellman equations and the associated discrete time stochastic optimal control problems.
\end{abstract}

\newpage

\tableofcontents

\section{Introduction}

\emph{Reinforcement learning} is an active research field in machine learning and has important applications in many areas which involve sequential decision making such as economics, engineering, finance, healthcare, logistics, and robotics (see, e.g., Sutton \& Barto \cite{sutton2018reinforcement}, Bertsekas \cite{bertsekas2019reinforcement}, and Bertsekas \& Tistsiklis \cite{bertsekas1996neuro} for overviews of the field and its application areas). In particular the combination of deep neural networks and reinforcement learning, i.e., deep reinforcement learning, 
has achieved remarkable success in complex decision-making problems during recent years
(see, e.g., Li \cite{li2017deep}, Fran{\c{c}}ois-Lavet et al. \cite{franccois2018introduction}, and Arulkumaran et al. \cite{arulkumaran2017deep} for survey articles).
The mathematical foundations of reinforcement learning are provided by the theory of \emph{stochastic optimal control} and, in particular, \emph{Markov decision processes} (MDPs; see, e.g., Bertsekas \& Shreve \cite{bertsekas1996stochastic}, Powell \cite{powell2007approximate}, and Puterman \cite{puterman2014markov}). 

The basis of many reinforcement learning algorithms such as temporal difference learning (see, e.g., Sutton \cite{sutton1988learning}), $Q$-learning (see, e.g., Watkins \cite{watkins1989learning} and Watkins \& Dayan \cite{watkins1992q}), and the SARSA algorithm (see, e.g., Rummery \& Niranjan \cite{rummery1994line}) is formed by the stochastic dynamic programming principle. 
It reduces the problem of making the optimal decision at a given state to solving a particular functional equation -- the so-called \emph{Bellman equation}. 
The approximative solution of Bellman equations in high dimensional state spaces
%Since these equations can almost never be solved explicitly their computational study constitutes an indispensable field of research. The design of efficient numerical schemes that perform well in high-dimensional state spaces 
is a notoriously difficult challenge due to the \emph{curse of dimensionality} (cf., e.g., Bellman \cite{Bellman57}, Novak \& Wozniakowski \cite[Chapter 1]{NovakWozniakowski2008I}, and Novak \& Ritter \cite{MR1485004}). 

In this work we introduce \emph{nonlinear Monte Carlo methods} for MDPs with infinite time horizon and finite control set that are polynomially tractable in the sense that the computational effort of the  algorithm to approximatively compute the solution of the Bellman equation grows at most polynomially in the reciprocal $1/\epsilon$ of the prescribed approximation accuracy $\epsilon \in (0,1]$ and the dimension $d\in \N=\{1,2,3,\ldots\}$ of the underlying state space. In particular, the proposed methods overcome the curse of dimensionality in the numerical approximation of solutions of Bellman equations. 

More formally, we consider MDPs that are specified by a measurable space $(\mathbb X,\mathcal X)$ (typically $(\R^d,\mathcal B(\R^d))$ with $d\in \N$ but our framework also covers discrete state spaces), a finite control set $A$, a discount factor $\delta \in (0,1)$, transition kernels $\kappa_a\colon \mathbb X \times \mathcal X \to [0,1]$, $a\in A$, and a measurable and bounded reward function $g\colon \mathbb X\times A \to \R$.
The set of strategies $\mathcal A$ consists of all measurable functions $\alpha\colon \mathbb X\to A$ (under appropriate assumptions the restriction to such non-randomized, stationary strategies is without loss of generality, see, e.g., \cite[Chapter 9]{bertsekas1996stochastic}). 
Every $\alpha \in \mathcal A$ defines a Markov process $X^\alpha=(X^\alpha_k)_{k\in \mathbb N_0}$ on a filtered probability space
$(\Omega, \mathcal F, (\mathcal F_k)_{k\in \mathbb N_0}, \mathbb{P})$ with state space $(\mathbb X,\mathcal X)$ such that $\mathbb{P}[X^\alpha_{k+1}\in B|\mathcal F_k]=\kappa_{\alpha(X^\alpha_k)}(B,X^\alpha_k)$ for all $k\in \N_0$, $B\in \mathcal X$. The expected gain of a strategy $\alpha \in \mathcal A$ when starting in $x\in \mathbb X$ is 
\begin{equation}
J(x,\alpha)=\mathbb E \!\left[\sum_{k=0}^\infty \delta^k g(X^\alpha_k,\alpha(X^\alpha_k))\Big |X^\alpha_0=x\right].
\end{equation}
The value function $v\colon \mathbb X \to \R$ is given by $v(x)=\sup_{\alpha \in \mathcal A}J(x,\alpha)$, $x\in \mathbb X$. The dynamic programming principle ensures that under appropriate conditions $v$ satisfies for all $x\in \mathbb X$ that
\begin{equation}\label{eq:Bellman_intro}
v(x)=\max_{a\in A}\left\{g(x,a)+\delta\, \mathbb E[v(X^{x,a})]\right\},
\end{equation}
where it holds for all $x\in \mathbb X$, $a\in A$ that the random variable $X^{x,a}\colon \Omega \to \mathbb X$ has distribution $\kappa_a(x,\cdot)$ (see, e.g., \cite[Chapter 9]{bertsekas1996stochastic} for a proof). 
An important subclass of these MDPs is given by optimal stopping problems in discrete time with infinite time horizon. For such problems the control set $A$ consists of the two elements ``stop'' and ``continue'' and once the control ``stop'' is taken the Markov process jumps to a hold state from which it cannot leave and generates no further gains.

There is a large number of numerical approximation methods for Bellman equations of the form \eqref{eq:Bellman_intro} that have been proposed and analyzed in the scientific literature.
%, see, e.g., Rust \cite{rust2008dynamic} for an overview. 
Deterministic numerical approximation methods for Bellman equations suffer in general from the curse of dimensionality. Indeed, Chow \& Tsitsiklis \cite{chow1989complexity,chow1991optimal} show for all $\dm \in \N$ that in the case where the state space $\mathbb X$ is given by the $d$-dimensional unit cube $[0,1]^d$ deterministic numerical approximation methods need at least $O(\epsilon^{-2d})$ 
computational operations to approximate the solution of the Bellman equation with precision $\epsilon\in (0,1]$. 
Rust \cite{rust1997using} shows that under certain assumptions it is possible to overcome this curse of dimensionality by allowing for randomized algorithms. Rust's method consists of randomly sampling grid points in the state space and performing value iteration on this stochastic grid. The method is model-based as it requires explicit knowledge of the transition density.
In the companion paper \cite{rust1997comparison} Rust points out that convergence of his method might break down in situations where the transition density has spikes. Also Kristensen et al.\ \cite{kristensen2021solving} report a dramatic rise in the method's variance as the state dimension increases. This is theoretically confirmed by 
Bray \cite{bray2022comment} who proves that Rust's method only overcomes the curse of dimensionality in the special case where the MDP is equivalent to an MDP where all but a vanishingly small fraction of state variables behave like history-independent uniform random variables. Variants of Rust's method for optimal stopping problems are designed by Broadie \& Glasserman in \cite{broadie1997pricing} and \cite{broadie2004stochastic}. The literature on optimal stopping comprises a variety of further randomized algorithms. We refer, for example, to 
%Longstaff \& Schwartz: 
\cite{longstaff2001valuing}, 
%Tsitsiklis \& van Roy: 
\cite{tsitsiklis1999optimal}, and \cite{tsitsiklis2001regression}
for regression-based algorithms for optimal stopping problems, we refer, for example, to
%Rogers: 
\cite{rogers2002monte},
%Anderson \& Broadie: 
\cite{andersen2004primal},
%Haugh \& Kogan
\cite{haugh2004pricing},
\cite{desai2012pathwise}, and
%\cite{belomestny2013multilevel}, 
\cite{belomestny2009true}
for duality-based algorithms for optimal stopping problems, 
and we refer, for example to \cite{becker2019deep}, \cite{becker2021solving}, and \cite{gonon2022deep} for deep learning-based algorithms for optimal stopping problems.

In this paper we introduce nonlinear Monte Carlo algorithms that overcome the curse of dimensionality in the approximation of Bellman equations, that require only very weak regularity assumptions and that are model-free in the sense that they do not use explicit knowledge of the transition kernel $\kappa$ but only need access to independent realizations of the one step transitions $X^{x,a}$ for arbitrary actions $a\in A$ and states $x\in \mathbb X$. Our approach combines ideas from \emph{$Q$-learning} and the recently introduced \emph{full-history recursive multilevel Picard} (MLP) approximations which have been proven to overcome the curse of dimensionality in the numerical approximation of certain semilinear partial differential equations (PDEs) (see, e.g., \cite{EHutzenthalerJentzenKruse2021, HJKNW2018, becker2020numerical, EHutzenthalerJentzenKruse2017, hutzenthaler2019overcoming, giles20019generalised, HJKN20, HutzenthalerKruse2017, hjk2019overcoming, beck2019overcoming, beck2020overcoming,beck2020nonlinear}). $Q$-learning is based on the idea to switch the order of expectation and maximization in \eqref{eq:Bellman_intro}. Formally, the $Q$-function satisfies for all $x\in \mathbb X$, $a\in A$ that $Q(x,a)=g(x,a)+\delta\, \E[v(X^{x,a})]$. This together with \eqref{eq:Bellman_intro} implies for all $x \in \X$ that $v(x)=\max_{a\in A}Q(x,a)$ and, hence, that for all $x \in \X$, $a \in A$ it holds that 
\begin{equation}\label{eq:dp_Q}
Q(x,a)=g(x,a)+\delta\, \E \!\left[\max_{b\in A}Q(X^{x,a},b)\right].
\end{equation}
Under appropriate conditions one can show that $Q$ is the unique solution of this fixed-point equation and that the sequence of fixed-point iterates $Q_n\colon \mathbb X \times A \to \R$, $n\in \N_0$, which satisfies for all $n\in \N$, $x\in \mathbb X$, $a\in A$ that $Q_0(x,a)= g(x,a)$ and
\begin{equation}\label{eq:Q_iterates} 
Q_{n}(x,a)=g(x,a)+\delta \, \E \!\left[\max_{b\in A}Q_{n-1}(X^{x,a},b)\right]
\end{equation}
converges to $Q$. We next employ a central idea of MLP approximations and decompose the iterates into multilevels to obtain for all $n\in \N_0$, $x\in \mathbb X$, $a\in A$ that
\begin{equation}
\begin{split}
Q_{n}(x,a)
%&=(TQ_{n-1})(a,x)=(TQ_0)(a,x)+\sum_{l=1}^{n-1}(TQ_{l})(a,x)-(TQ_{l-1})(a,x)\\
&=g(x,a)+\delta \left[ \sum_{l=0}^{n-1} \E \!\left[\!\left(\!\max_{b\in A}Q_{l}(X^{x,a},b)\!\right)- \mathbbm{1}_{\N}(l)\! \left(\!\max_{b\in A}Q_{l-1}(X^{x,a},b)\!\right)\! \right]\right].
\end{split}
\end{equation}
In this telescope expansion, we apply a fundamental idea of Heinrich \cite{heinrich1998monte,heinrich2001multilevel} and Giles \cite{giles2008improved} and approximate the expected values by Monte Carlo averages with different degrees of accuracy at different levels $l\in \{1,\ldots,n\}$. The convergence of $(Q_l)_{l\in \N}$ ensures that for large $l\in \{1,\ldots,n\}$ the difference between $Q_{l}$ and $Q_{l-1}$ is small and hence we use for large $l\in \{1,\ldots,n\}$ less Monte Carlo samples to approximate the expected value $\E\left[\left(\max_{b\in A}Q_{l}(X^{x,a},b)\right)-\left(\max_{b\in A}Q_{l-1}(X^{x,a},b)\right)\right]$ than for small $l\in \{1,\ldots,n\}$. More specifically, we fix $M\in \N$ and use $M^{n-l}$ independent Monte Carlo samples to approximate the expected value at level $l\in \{1,\ldots,n\}$ which leads to the full-history recursive multilevel fixed-point (MLFP) approximation scheme in \eqref{eq:scheme_intro} below.

To briefly sketch the contribution of this article within this introductory section, we now present in the following result, \cref{mlp_approx_cor_ctrl_slim} below, a special case of \cref{mlp_approx_main_theo}, the main result of this article. Below  \cref{mlp_approx_cor_ctrl_slim} we explain in words the statement of  \cref{mlp_approx_cor_ctrl_slim} as well as the mathematical objects appearing in \cref{mlp_approx_cor_ctrl_slim}.

\begin{theo} \label{mlp_approx_cor_ctrl_slim} 
Let $A$ be a nonempty set, for every $\dm \in \N$ let $g_\dm \colon \R^d \times A \rightarrow \R$ be $(\mathcal{B}(\R^d) \otimes \pset{A}) / \mathcal{B}(\R)$-measurable, assume $\sup_{\dm \in \N} \sup_{x \in \R^\dm}\sup_{ a\in A} |g_\dm(x, a)| < \infty$, let $(\Omega, \mathcal{F}, \mathbb{P})$ be a probability space, let $\Theta = \cup_{n \in \N} \Z^n$, for every $\dm \in \N$ let $\mathcal{F}_\dm^\theta \subseteq \mathcal{F}$, $\theta \in \Theta$, be independent sub-sigma-algebras of $\mathcal{F}$, for every $\dm \in \N$ let $X_d^{\theta} =  (X_\dm^{\theta, x, a} (\omega) )_{ (x, a, \omega) \in \R^\dm \times A\times \Omega } \colon \R^\dm \times A \times \Omega \rightarrow \R^\dm$, $\theta \in \Theta,$ be i.i.d.\ random fields which satisfy for all $\dm \in \N$, $\theta \in \Theta$ that $X_\dm^\theta$ is $ ( \mathcal{B}(\R^\dm) \otimes \pset{A} \otimes \mathcal{F}_\dm^\theta ) /  \mathcal{B}(\R^\dm)$-measurable, for every $\dm \in \N$ let $\delta_\dm \in [0,1)$, $\mathcal{R}_\dm\in [0, \infty)$,
let $M \in \N \cap [\sup_{\dm \in \N} (4 |A|^2 (1-\delta_\dm)^{-2}), \infty]$,
for every $\dm \in \N$ let $\mathcal{Q}_{\dm, n}^\theta \colon \R^\dm \times A \times \Omega \rightarrow \R$, $n \in \N_0$, $\theta \in \Theta$, satisfy for all $n \in \N_0$, $\theta \in \Theta$, $x \in \R^\dm$, $a\in A$ that
\begin{multline}\label{eq:scheme_intro}
\mathcal{Q}_{\dm, n}^{\theta}(x,a) = g_\dm(x,a) \\+ \sum_{l = 0}^{n-1} \frac{\delta_\dm}{M^{n-l}} \sum_{i = 1}^{M^{n-l}} \left[ \max_{b \in A}\big\{\! \mathcal{Q}_{\dm, l}^{(\theta, l, i)}(X_\dm^{(\theta, l, i), x,a},b) \! \big\} -\mathbbm{1}_{\N}(l) \max_{b \in A} \big\{\! \mathcal{Q}_{\dm, \max\{ l-1, 0 \}}^{(\theta, -l, i)}(X_\dm^{(\theta, l, i), x,a},b) \!\big\}\right], 
\end{multline}
and for every $\dm \in \N$ let $\mathcal{C}_{\dm, n} \in [0, \infty)$, $n \in \N_0$, satisfy for all $n \in \N_0$ that
\begin{align} \label{eq:comp_cost_intro}
\mathcal{C}_{\dm, n} \leq \sum_{l = 0}^{n-1} M^{n-l}\big( \mathcal{R}_\dm + \mathcal{C}_{\dm, l} + \mathcal{C}_{\dm, \max \{l-1, 0\}} \mathbbm{1}_{\N}(l)\big).
\end{align}
Then
\begin{enumerate}
\item[\textnormal{(i)}]
\label{it:ex_uni_bellman_intro} it holds for all $\dm \in \N$ that there exists a unique bounded $(\mathcal{B}(\R^\dm)\otimes \pset{A})/ \mathcal{B}(\R)$-measurable $Q_\dm \colon \R^\dm \times A \rightarrow \R$ which satisfies for all $x \in \R^\dm$, $a \in A$ that 
\begin{align}
Q_\dm(x,a) = g_\dm(x, a) + \delta_\dm\, \E\! \left[ \max_{b \in A} Q_d(X_\dm^{0,x,a}, b) \right]
\end{align}
and
\item[\textnormal{(ii)}] there exist $N \colon (0, 1] \rightarrow \N$ and $c \in \R$ such that for all $\dm \in \N$, $\varepsilon \in (0,1]$ it holds that $\mathcal{C}_{\dm, N_\varepsilon} \leq c \mathcal{R}_\dm \varepsilon^{-c}$ and
\begin{align}
\sup_{x \in \R^\dm} \left( \E\!\left[ \max_{a \in A}| Q_\dm(x,a) - \mathcal{Q}_{\dm, N_\varepsilon}^0(x, a) |^2 \right] \right)^{\!\!\nicefrac{1}{2}} \leq \varepsilon.
\end{align}
\end{enumerate}
\end{theo}

\cref{mlp_approx_cor_ctrl_slim} is an immediate consequence from \cref{mlp_approx_cor_ctrl} in \cref{sec:comp_comp} below. \cref{mlp_approx_cor_ctrl}, in turn, follows from \cref{mlp_approx_main_theo}, which is the main result of this article. In the following we add some comments on the mathematical objects appearing in \cref{mlp_approx_cor_ctrl_slim} above.

In \cref{mlp_approx_cor_ctrl_slim} we introduce in \eqref{eq:scheme_intro} a Monte Carlo-type approximation algorithm for a sequence of MDPs indexed by the dimension $d\in \N$ of the state space. To formulate the proposed Monte Carlo-type approximation algorithm in \eqref{eq:scheme_intro} we need, roughly speaking, sufficiently many independent random variables which are indexed over a sufficiently large index set. This sufficiently large index set is provided through the set $\Theta = \cup_{n \in \N} \Z^n$ introduced
%in the first line of 
in \cref{mlp_approx_cor_ctrl_slim}. The triple $(\Omega, \mathcal{F}, \mathbb{P})$ in 
%
%the first line of 
%
\cref{mlp_approx_cor_ctrl_slim} is the probability space on which the random variables are defined. In \cref{mlp_approx_cor_ctrl_slim} we consider for every $d\in \N$ an MDP with state space $(\R^d, \mathcal B(\R^d))$. We assume that all elements of the sequence of MDPs have a common control set $A$.  
For every $d\in \N$, $x\in \R^d$, $a\in A$ the one-step transition of the controlled Markov chain of the MDP with state space $(\R^d, \mathcal B(\R^d))$ is given by the random variable $X_d^{0,x,a}\colon \Omega \to \R^d$. 
%
%
%that is introduced in the fourth line of \cref{mlp_approx_cor_ctrl_slim}. 
%
%
In the language of MDPs for every $d\in \N$, $a\in A$ the transition kernel $\kappa_{d,a}\colon \R^d\times \mathcal B(\R^d) \to [0,1]$ is thus determined by the distribution of $X_d^{0,\cdot,a}$. For every $d\in \N$ the function $g_d\colon \R^d\times A \to \R$ introduced in the first line of \cref{mlp_approx_cor_ctrl_slim} is the reward function of the MDP with state space $(\R^d, \mathcal B(\R^d))$. In \cref{mlp_approx_cor_ctrl_slim} we assume that the functions $g_d$, $d\in \N$, are uniformly bounded in $d\in \N$. For every $d\in \N$ the real number $\delta_d\in [0, 1)$ introduced in
%
%
%the sixth line of
%
%
\cref{mlp_approx_cor_ctrl_slim} is the discount factor of the MDP with state space $(\R^d, \mathcal B(\R^d))$. 

Item (i) in \cref{mlp_approx_cor_ctrl_slim} establishes the essentially well-known result that under the above assumptions for every $d\in \N$ the Bellman equation \eqref{eq:Bellman_intro} for the $Q$-function associated to the MDP with state space $(\R^d, \mathcal B(\R^d))$ has a unique solution $Q_d\colon \R^d\times A\to \R$.

In \eqref{eq:scheme_intro} in \cref{mlp_approx_cor_ctrl_slim} we specify the MLFP approximation scheme which we propose to approximate the solution of the Bellman equation \eqref{eq:Bellman_intro}. The MLFP approximations $\mathcal{Q}^\theta_{\dm, n}$, $\dm\in \N$, $n\in \N_0$, $\theta \in \Theta$, are indexed by the dimension $d\in \N$, by the number $n\in \N_0$ of fixed-point iterates and by a parameter $\theta \in \Theta$ which is different for different appearences of MLFP approximations in \eqref{eq:scheme_intro}. As random input sources the MLFP approximation scheme proposed in \eqref{eq:scheme_intro} employs the random variables $X_d^{\theta,x,a}\colon \Omega \to \R^d$, $\theta \in \Theta \setminus \{0\}$, $d\in \N$, $x\in \R^d$, $a\in A$. Note that for every $\theta \in \Theta \setminus \{0\}$, $d\in \N$, $x\in \R^d$, $a\in A$ the random variable $X_d^{\theta,x,a}\colon \Omega \to \R^d$ which is used as random input source of the MLFP approximation scheme proposed in \eqref{eq:scheme_intro} and the random variable $X_d^{0,x,a}\colon \Omega \to \R^d$ which is used to formulate the Bellman equation \eqref{eq:Bellman_intro} are identically distributed. The parameter $\theta \in \Theta$ ensures that different appearances of MLFP approximations in \eqref{eq:scheme_intro} are independent and this ensures that \eqref{eq:scheme_intro} can be implemented with recursive function calls.
The natural number $M \in \N \cap [\sup_{\dm \in \N} (4 |A|^2 (1-\delta_\dm)^{-2}), \infty] $ in \cref{mlp_approx_cor_ctrl_slim} determines the number of Monte Carlo samples used in the definition of the MLFP approximation $\mathcal{Q}^\theta_{\dm,n}$, $d\in \N$, $n\in \N_0$, $\theta \in \Theta$ in \eqref{eq:scheme_intro}. The assumption that the natural number $M$ is an element of the set $ \N \cap [\sup_{\dm \in \N} (4 |A|^2 (1-\delta_\dm)^{-2}), \infty] $ ensures that $ \sup_{\dm \in \N} (4 |A|^2 (1-\delta_\dm)^{-2}) < \infty$. This implies that the control set $A$ is finite and that the sequence $(\delta_\dm)_{\dm\in \N}\subseteq [0,1)$ is bounded away from $1$ in the sense that $\sup_{\dm \in \N} \delta_\dm < 1$.
Note that 
by a suitable identification 
this framework also covers the case where for each $d\in \N$ there is an individual control set $A_d$ such that $\sup_{d\in \N}|A_d|<\infty$.

For every $d\in \N$ the nonnegative real number $\mathcal{R}_\dm \in [0, \infty)$ in \cref{mlp_approx_cor_ctrl_slim} is understood as an upper bound of the computational cost to compute one realization of any of the random variables $X^{\theta,x}_d\colon \Omega \to (\R^d)^A$, $\theta \in \Theta$, $x\in \R^d$. The real numbers $\mathcal{C}_{\dm,n}\in [0,\infty)$, $\dm \in \N$, $n\in \N_0$, in \eqref{eq:comp_cost_intro} in \cref{mlp_approx_cor_ctrl_slim} model the computational costs of the MLFP approximation scheme in \eqref{eq:scheme_intro}. More specifically, 
for every $d\in \N$, $n\in \N_0$ the real number $\mathcal{C}_{\dm, n}\in [0,\infty)$ 
%in \eqref{eq:comp_cost_intro} in \cref{mlp_approx_cor_ctrl_slim}
represents an upper bound of the computational costs to compute the realizations of all random variables $X^{\theta,x}_d\colon \Omega \to (\R^d)^A$, $\theta \in \Theta$, $x\in \R^d$, required to compute one realization of $\mathcal{Q}^0_{\dm, n}(0,\cdot)\colon \Omega \to \R^A$.

Item (ii) in \cref{mlp_approx_cor_ctrl_slim} proves that the solutions of the Bellman equations in \eqref{eq:Bellman_intro} can be approximated by means of the MLFP approximation scheme in \eqref{eq:scheme_intro} with a computational cost which grows at most polynomially in the reciprocal $1/\epsilon$ of the prescribed approximation accuracy $\epsilon\in (0,1]$ and linearly in the computational cost $\mathcal{R}_d$ to compute one realization of any of the random variables $X^{\theta,x}_d\colon \Omega \to (\R^d)^A$, $\theta \in \Theta$, $x\in \R^d$, where $d\in \N$ is the dimension of the state space of the associated MDP. In particular, if the computational cost to compute one realization of any of the random variables $X^{\theta,x}_d\colon \Omega \to (\R^d)^A$, $\theta \in \Theta$, $x\in \R^d$, grows at most polynomially in the dimension $d\in \N$ of the state space of the associated MDP (as it is often the case in practical applications), then 
the MLFP approximation scheme in \eqref{eq:scheme_intro} overcomes the curse of dimensionality for the approximation of the solutions of the Bellman equations in \eqref{eq:Bellman_intro}. However, we would like to point out that the constant $c \in \R$ appearing in item~(ii) in \cref{mlp_approx_cor_ctrl_slim} may become arbitrary large if $\sup_{\dm \in \N}\delta_\dm$ is close to 1. Thus the computational cost of the nonlinear Monte Carlo methods in \eqref{eq:scheme_intro} may become impractical even so the methods in \eqref{eq:scheme_intro} provably overcome the curse of dimensionality.

In the following we also add some comments on generalizations and variants of \cref{mlp_approx_cor_ctrl_slim} presented in this article. 
 While \cref{mlp_approx_cor_ctrl_slim} considers a sequence of MDPs indexed by the dimension $d\in \N$ of the Euclidean state spaces $(\R^d, \mathcal B(\R^d))$, \cref{mlp_approx_cor_ctrl} considers a family of MDPs with general index set $\Dm$ and general measurable state spaces $(\mathbb X_d, \mathcal X_d)$, $d\in \Dm$.  
While \cref{mlp_approx_cor_ctrl_slim} considers MDPs with bounded reward functions, \cref{mlp_approx_cor_ctrl} allows for unbounded reward functions. 
While \cref{mlp_approx_cor_ctrl_slim} and \cref{mlp_approx_cor_ctrl} consider Bellman equations of MDPs, \cref{mlp_approx_main_theo} considers more general functional fixed-point equations (we refer to \eqref{eq:gen_func_equation} in  \cref{mlp_approx_main_theo} for details). \Cref{mlp_approx_cor_stop_slim} proves that a variant of the MLFP approximation scheme (see \eqref{mlp_approx_cor_stop_slim_scheme} for details) overcomes the curse of dimensionality for the approximation of the solutions of Bellman equations for optimal stopping problems (see \eqref{eq:Bellman_stopping} for details).

The remainder of this article is organized as follows. In \cref{sec:ex_uniq_sol_fp_eq} below we establish existence, uniqueness, and integrability properties for solutions of functional fixed-point equations. In \cref{sec:mlfp_approx} below we introduce MLFP approximations for solutions of functional fixed-point equations, we study measurability, distributional, and integrability properties for the introduced MLFP approximations and we establish recursive and subsequently non-recursive upper bounds for the $L^2$-distances between the exact solutions of the considered functional fixed-point equations and the proposed MLFP approximations. In \cref{sec:comp_comp} we combine the existence, uniqueness, and regularity properties for solutions of functional fixed-point equations, which we have established in \cref{sec:ex_uniq_sol_fp_eq}, with the error analysis for MLFP approximations for functional fixed-point equations, which we have established in \cref{sec:mlfp_approx}, to obtain a computational complexity analysis for MLFP approximations for functional fixed-point equations and for Bellman equations of MDPs and optimal stopping problems.

\section{Existence and uniqueness of solutions of functional fixed-point equations} \label{sec:ex_uniq_sol_fp_eq}

\begin{defi} \label{defi_fam_kernel}
Let $(\X, \mathcal{X})$, $(\Y, \mathcal{Y})$ be a nonempty measurable spaces, let $A$ be a nonempty set, and let $\kappa_a \colon \X \times \mathcal{Y} \rightarrow [0,1]$, $a \in A$, satisfy for all $a \in A$, $M \in \mathcal{Y}$ that $\X \ni x \mapsto \kappa_a(x,M) \in [0,1]$ is $\mathcal{X}/\mathcal{B}\big( [0,1] \big)$-measurable and for all $x \in \X$, $a \in A$ that $\mathcal{Y} \ni M \mapsto \kappa_a(x,M) \in [0,1]$ is a probability measure on $(\Y, \mathcal{Y})$. Then we say that $(\kappa_a)_{a \in A}$ is a family of stochastic kernels from $(\X, \mathcal{X})$ to $(\Y, \mathcal{Y})$.
\end{defi}

\begin{lemma} \label{ex_unique_sol_sto_fp_eq_ker} Let $c,L \in [0,\infty)$ with $cL < 1$, let $(\mathbb{X}, \mathcal{X})$ be a nonempty measurable space, let $A$ be a nonempty countable set, let $ \kappa = (\kappa_a)_{a\in A}$ be a family of stochastic kernels from $(\X, \mathcal{X})$ to $(\X, \mathcal{X})$, let $\R^A = \{ r \colon A \rightarrow \R \}$, let $\A = \bigotimes_{a \in A} \mathcal{B}(\R)$, let $f \colon \X \times \R^A \rightarrow \R$ be $(\mathcal{X} \otimes \A)/\mathcal{B}(\R)$-measurable, let $\w \colon \X \rightarrow (0,\infty)^A$ be $\mathcal{X}/\bigotimes_{a \in A}\mathcal{B}\big((0, \infty)\big)$-measurable, let
\begin{align}
\W = \bigg\{\! (u \colon \X \rightarrow \R^A) : u \text{ is } \mathcal{X}/\A\text{-measurable, } \sup_{(x,a) \in \X \times A} \big[|\dkl{\w}{x}{a}|^{-1} |\dkl{u}{x}{a}| \big] < \infty \! \bigg\},
\end{align}
assume for all $x \in \X$, $a \in A$, $r,s \in \R^A$ that $|f(x,r) - f(x,s)| \leq L \sup_{b \in A} |r(b) - s(b) |$, $\int_\X \sup_{b \in A} \dkl{\w}{y}{b}\kappa_a(x,dy) \leq c \dkl{\w}{x}{a},$ and $\sup_{(t,b) \in \X \times A} [|\dkl{\w}{t}{b}|^{-1} \int_\X |f(y,0)|\kappa_b(t,dy)] < \infty$. Then there exists a unique $v \in \W$ which satisfies for all $x \in \X$, $a \in A$ that 
\begin{align}
\int_\X |f(y,v(y))| \kappa_a(x,dy) < \infty\qquad \text{and} \qquad \dkl{v}{x}{a} = \int_\X f(y,v(y))\kappa_a(x,dy).
\end{align}
\end{lemma}

\begin{mproof}{\cref{ex_unique_sol_sto_fp_eq_ker}}
Let $\left\| \cdot \right\| \colon \W \rightarrow [0,\infty)$ satisfy for all $u \in \W$ that $ \|u\| =  \sup_{(x,a) \in \X \times A} \frac{|\left(u(x)\right)(a)|}{|\left(\w (x)\right)(a)|} $. Note that $\left\| \cdot \right\|$ is a norm on $\W$ and $(\W, \left\| \cdot \right\|)$ is a Banach space. The assumption that for all $x \in X$, $r,s \in \R^A$ it holds that $|f(x,r) - f(x,s)| \leq L \sup_{a \in A} |r(a) - s(a)|$ and the assumption that for all $x \in \X$, $a \in A$ it holds that $\int_\X \sup_{b \in A} \dkl{\w}{y}{b} \kappa_a(x,dy) \leq c \dkl{\w}{x}{a} $, yield that for all $u \in \W$, $x \in \X$, $a \in A$ it holds that
\begin{align}
&\hspace{-0.5cm}\frac{1}{\dkl{\w}{x}{a} } \bigg| \int_\X f(y,u(y)) \kappa_a(x,dy) \bigg| \nonumber \\
&\leq \frac{1}{\dkl{\w}{x}{a}} \int_\X |f(y,u(y)) - f(y,0)| \kappa_a(x,dy) + \frac{1}{\dkl{\w}{x}{a}} \int_\X |f(y,0)| \kappa_a(x,dy)\nonumber\\
&\leq \frac{L}{\dkl{\w}{x}{a}} \int_\X \sup_{b \in A} |\dkl{u}{y}{b}| \kappa_a(x,dy) + \frac{1}{\dkl{\w}{x}{a}} \int_\X |f(y,0)| \kappa_a(x,dy)\nonumber\\
&= \frac{L}{\dkl{\w}{x}{a}} \int_\X \sup_{b \in A} \bigg\{ \frac{| \dkl{u}{y}{b}|}{\dkl{\w}{y}{b}}\dkl{\w}{y}{b} \bigg\} \kappa_a(x,dy) + \frac{1}{\dkl{\w}{x}{a}} \int_\X |f(y,0)| \kappa_a(x,dy)\nonumber\\
&\leq \frac{L \| u \|}{\dkl{\w}{x}{a}} \int_\X \sup_{b \in A}\dkl{\w}{y}{b} \kappa_a(x,dy) + \frac{1}{\dkl{\w}{x}{a}} \int_\X |f(y,0)| \kappa_a(x,dy)\nonumber\\
&\leq cL \|u\| + \frac{1}{\dkl{\w}{x}{a}} \int_\X |f(y,0)| \kappa_a(x,dy).
\end{align}
The assumption that $\sup_{(x,a) \in \X \times A} \big[|\dkl{\w}{x}{a}|^{-1}\int_\X \big| f(x,0) \big| \kappa_a(x,dy)\big] < \infty$ demonstrates that for all $u \in \W$ it holds that
\begin{align}
\sup_{(x,a) \in \X \times A} \frac{1}{\dkl{\w}{x}{a} } \bigg| &\int_\X f(y,u(y)) \kappa_a(x,dy) \bigg|\nonumber\\
&\leq cL \|u \| + \sup_{(x,a) \in \X \times A} \frac{1}{\dkl{\w}{x}{a}} \int_\X |f(y,0)| \kappa_a(x,dy) < \infty.
\end{align}
\cite[Lemma 14.20]{KlenkeProb2014} ensures that for all $u \in \W$ it holds that the map $\X \ni x \mapsto \big[ A \ni a \mapsto \int_\X f(y,u(y)) \kappa_a(x,dy) \in \R \big] \in \R^A$ is $\mathcal{X}/\A$-measurable. 
Let $\Phi \colon \W \rightarrow \W$ be the function which satisfies for all $u \in \W$, $(x,a) \in \X \times A$ that
\begin{align}
[\Phi(u)](x)(a) = \int_\X f(y,u(y))\kappa_a(x,dy).
\end{align}
The assumption that for all $x \in \X$, $r,s \in \R^A$ it holds that $|f(x,r) - f(x,s)| \leq L \sup_{a \in A} |r(a) - s(a)|$, the assumption that $\sup_{(x,a) \in \X \times A} \big[|\dkl{\w}{x}{a}|^{-1}\int_\X \big| f(x,0) \big| \kappa_a(x,dy)\big] < \infty$, and the assumption that for all $(x,a) \in \X \times A$ it holds that $\int_\X \sup_{b \in A}\dkl{\w}{y}{b} \kappa_a(x,dy) \leq c \dkl{\w}{x}{a}$, ensure that for all $u,v \in \W, (x,a) \in \X \times A$ it holds that
\begin{align}
\frac{1}{ \dkl{\w}{x}{a}}\Big| [\Phi(u)](x)(a) &- [\Phi(v)](x)(a) \Big|\nonumber\\
&\leq \frac{1}{\dkl{\w}{x}{a}} \int_\X \big| f(y,u(y)) - f(y,v(y)) \big|\kappa_a(x,dy)\nonumber\\
&\leq \frac{L}{\dkl{\w}{x}{a}} \int_\X \sup_{b \in A} | \dkl{u}{y}{b} - \dkl{v}{y}{b} | \kappa_a(x,dy)\nonumber\\
&= \frac{L}{\dkl{\w}{x}{a}} \int_\X \sup_{b \in A} \bigg\{ \frac{|\dkl{u}{y}{b} - \dkl{v}{y}{b}|}{\dkl{\w}{y}{b}} \dkl{\w}{y}{b} \bigg\} \kappa_a(x,dy)\nonumber\\
&\leq \frac{L \| u - v \|}{\dkl{\w}{x}{a}} \int_\X \sup_{b \in A} \dkl{\w}{y}{b}\kappa_a(x,dy)\nonumber\\
&\leq cL \| u - v \|.
\end{align}
This in turn proves that for all $u,v \in \W$ it holds that
\begin{align}
\big\| \Phi(u) - \Phi(v) \big\| \leq cL\|u - v \|.
\end{align}
The assumption that $cL < 1$ shows that $\Phi$ is a contraction. Hence Banach's fixed-point theorem proves that there exists a unique function $v \in \W$ such that $v = \Phi(v)$. The fact that $v \in \W$ ensures that for all $(x,a) \in \X \times A$ it holds that $\int_\X | f(y,v(y)) | \kappa_a(x,dy) < \infty$. The proof of \cref{ex_unique_sol_sto_fp_eq_ker} is thus completed.
\end{mproof}

\begin{corollary} \label{ex_unique_sol_sto_fp_eq_ew} Let $c,L \in [0, \infty)$ with $cL < 1$, let $(\X, \mathcal{X})$ be a nonempty measurable space, let $(\Omega, \mathcal{F}, \mathbb{P})$ be a probability space, let $A$ be a nonempty countable set, let $\R^A = \{ r \colon A \rightarrow \R \}$, let $\A = \bigotimes_{a \in A} \mathcal{B}(\R)$, let $\X^A = \{ q \colon A \rightarrow \X \}$, let $\wt{\mathcal{X}} = \bigotimes_{a \in A} \mathcal{X}$, let $X  = \big( X^{x,a}(\omega) \big)_{x \in \X,\;a\in A,\; \omega \in \Omega}\colon \X \times \Omega \rightarrow \X^A$ be $(\mathcal{X} \otimes \mathcal{F})/\wt{\mathcal{X}}$-measurable, let $f \colon \X \times \R^A \rightarrow \R$ be $(\mathcal{X} \otimes \A)/\mathcal{B}(\R)$-measurable, let $\w \colon \X \rightarrow (0,\infty)^A$ be $\mathcal{X}/\bigotimes_{a \in A}\mathcal{B}\big((0, \infty)\big)$-measurable, let 
\begin{align*}
\W = \big\{ (u \colon \X \rightarrow \R^A) : u \text{ is } \mathcal{X}/\A\text{-measurable, } \sup_{(x,a) \in \X \times A} \big[|\dkl{\w}{x}{a}|^{-1} |\dkl{u}{x}{a}| \big] < \infty \big\},
\end{align*}
assume that for all $(x,a) \in \X \times A$, $r,s \in \R^A$ it holds that $|f(x,r) - f(x,s)| \leq L \sup_{b \in A} |r(b) - s(b)|$, $\E [\sup_{b \in A} \dkl{\w}{X^{x,a}}{b}] \leq c \dkl{\w}{x}{a}$, and $\sup_{(y, b) \in \X \times A} \big[ | \dkl{\w}{y}{b}|^{-1}\E[|f(X^{y,b},0)|] \big] < \infty$. Then there exists a unique function $v \in \W$ which satisfies for all $(x,a) \in \X \times A$ that
\begin{align}
\E \big[|f(X^{x,a},v(X^{x,a}))|\big] < \infty, \quad \text{and} \quad \dkl{v}{x}{a} = \E \big[ f(X^{x,a},v(X^{x,a})) \big].
\end{align}
\end{corollary}

\begin{mproof}{\cref{ex_unique_sol_sto_fp_eq_ew}}
Let $\kappa_a \colon \X \times \mathcal{X} \rightarrow [0,1]$, $a \in A$, satisfy for all $(x,a) \in \X \times A$, $M \in \mathcal{X}$, that $\kappa_a(x,M) = \mathbb{P}[X^{x,a} \in M]$. The fact that $X$ is $(\mathcal{X} \otimes \mathcal{F})/\wt{\mathcal{X}}$-measurable ensures that for all $(x,a) \in \X \times A$ it holds that $X^{x,a} \colon \Omega \rightarrow \X$ is $\mathcal{F}/\mathcal{X}$-measurable. This implies that for all $(x,a) \in \X \times A$ it holds that $\mathcal{X} \ni M \mapsto \kappa_a(x,M) = \mathbb{P}[X^{x,a} \in M] \in [0,1]$ is a probability measure on $(\X ,\mathcal{X})$. The fact that $X$ is $(\mathcal{X} \otimes \mathcal{F})/\wt{\mathcal{X}}$-measurable and \cite[Theorem 14.16]{KlenkeProb2014} imply that for all $a \in A,$ $M \in \mathcal{X}$ it holds that $\X \ni x \mapsto \kappa_a (x,M) = \mathbb{P}[X^{x,a} \in M]$ is $\mathcal{X}/ \mathcal{B}([0,1])$-measurable. Hence it holds that $(\kappa_a)_{a \in A}$ is a familiy of stochastic kernels form $(\X, \mathcal{X})$ to $(\X, \mathcal{X})$. The assumptions that for all $(x,a)\in \X \times A$ it holds that $\E \big[ \sup_{b \in A} \dkl{\w}{X^{x,a}}{b} \big] \leq c \dkl{\w}{x}{a}$, and $\sup_{(y,b) \in \X \times A} \big[ | \dkl{\w}{y}{b} |^{-1} \E[ |f(X^{y,b}, 0 )| ] \big] < \infty$, prove that for all $(x,a) \in \X \times A$ it holds that 
\begin{align}
\int_\X \sup_{b \in A} \dkl{\w}{y}{b} \kappa_a(x,dy) &= \hspace{-0.15cm} \int_\X \sup_{b \in A} \dkl{\w}{y}{b} \big( X^{x,a}(\mathbb{P}) \big)(dy) = \E [\sup_{b \in A} \dkl{\w}{X^{x,a}}{b} ]\leq c \dkl{\w}{x}{a},
\end{align}
and
\begin{align}
\sup_{(y,b)\in \X \times A} \frac{1}{\dkl{\w}{y}{b}}  \int_\X |f(t,0)| \kappa_b(y,dt) &= \sup_{(y,b)\in \X \times A} \frac{1}{\dkl{\w}{y}{b}} \int_\X |f(t,0)|\big( X^{y,b}(\mathbb{P}) \big)(dt) \nonumber \\
&= \sup_{(y,b)\in \X \times A} \frac{1}{\dkl{\w}{y}{b}} \E \big[ |f(X^{y,b},0)| \big] < \infty.
\end{align}
\cref{ex_unique_sol_sto_fp_eq_ker} implies that there exists a unique function $v \in \W$ such that for all $(x,a) \in \X \times A$ it holds that
\begin{align}
\E\Big[ \big|f\big(X^{x,a}, v(X^{x,a})\big)\big| \Big] &= \int_\X | f( y, v(y) ) | \big( X^{x,a}(\mathbb{P}) \big)(dy) = \int_\X | f(y,v(y)) | \kappa_a(x,dy) < \infty,
\end{align}
and
\begin{align}
\dkl{v}{x}{a} &= \int_\X f(y,v(y)) \kappa_a(x,dy) = \int_\X f(y, v(y)) \big( X^{x,a}(\mathbb{P}) \big)(dy) = \E \big[ f\big(X^{x,a}, v(X^{x,a}) \big) \big].
\end{align}
The proof of \cref{ex_unique_sol_sto_fp_eq_ew} is thus completed.
\end{mproof}

\begin{lemma}\label{sol_sto_fp_eq_estimates} Let $c_f, c_\w, L \in [0,\infty)$ with $c_\w L < 1$, let $(\X, \mathcal{X})$ be a nonempty measurable space, let $(\Omega, \mathcal{F}, \mathbb{P})$ be a probability space, let $A$ be a nonempty countable set, let $\R^A = \{ r\colon A \rightarrow \R  \}$, let $\A = \bigotimes_{a \in A} \mathcal{B}(\R)$, let $\X^A = \{ q \colon A \rightarrow \X \}$, let $\wt{\mathcal{X}} = \bigotimes_{a \in A}\mathcal{X}$, let $X = \big( X^{x,a}(\omega)\big)_{x \in \X,\; a\in A,\; \omega \in \Omega} \colon \X \times \Omega \rightarrow \X^A$ be $(\mathcal{X} \otimes \mathcal{F})/\wt{\mathcal{X}}$-measurable, let $f \colon \X \times \R^A \rightarrow \R$ be $(\mathcal{X}\otimes \A)/\mathcal{B}(\R)$-measurable, let $\w \colon \X \rightarrow (0,\infty)^A$ be $\mathcal{X}/\bigotimes_{a \in A}\mathcal{B}((0,\infty))$-measurable, let $v \colon \X \rightarrow \R^A$ be $\mathcal{X}/\A$-measurable, assume that for all $(x,a)\in \X \times A$, $r,s \in \R^A$ it holds that $|f(x,r) - f(x,s)| \leq L \sup_{b \in A} |r(b) - s(b)|$, $\sup_{(y,b) \in \X \times A} \big[ | \dkl{\w}{y}{b}|^{-1} |\dkl{v}{y}{b} | \big] < \infty$, \hspace{0.0cm} $\big( \E\big[\sup_{b \in A} \big|\dkl{\w}{X^{x,a}}{b}\big|^2\big] \big)^\frac{1}{2} \leq c_\w \dkl{\w}{x}{a}$, $\big( \E\big[ | f(X^{x,a},0) |^2 \big] \big)^\frac{1}{2}\leq c_f \dkl{\w}{x}{a}$,  and $\dkl{v}{x}{a} = \E[f(X^{x,a},v(X^{x,a}))]$. Then it holds that
\begin{align}
\sup_{(x,a) \in \X \times A} \frac{| \dkl{v}{x}{a}|}{\dkl{\w}{x}{a}} \leq \frac{c_f}{1-c_\w L} \quad \text{and} \hspace{0.2cm} \sup_{(x,a) \in \X \times A} \frac{\big( \E\big[ \sup_{b \in A} | \dkl{v}{X^{x,a}}{b}|^2 \big] \big)^\frac{1}{2}}{\dkl{\w}{x}{a}} \leq \frac{c_f c_\w}{1- c_\w L}.
\end{align}
\end{lemma}

\begin{mproof}{\cref{sol_sto_fp_eq_estimates}}
Jensen's inequality, the triangle inequality, and the assumption that for all $x \in \X$, $r,s \in \R^A$ it holds that $|f(x,r) - f(x,s)| \leq L \sup_{a \in A} |r(a) - s(a)|$ prove that for all $(x,a) \in \X \times A$ it holds that
\begin{align}
\frac{|\dkl{v}{x}{a}|}{\dkl{\w}{x}{a}} &= \frac{\Big| \E \big[ f(X^{x,a}, v(X^{x,a})) \big] \Big|}{\dkl{\w}{x}{a}} \leq \frac{\E \big[ |f(X^{x,a}, v(X^{x,a}))|^2 \big]^\frac{1}{2}}{\dkl{\w}{x}{a}} \nonumber \\
&\leq L\frac{\E\big[ \sup_{b \in A}| \dkl{v}{X^{x,a}}{b} |^2 \big]^\frac{1}{2}}{\dkl{\w}{x}{a}} + \frac{\E \big[ |f(X^{x,a}, 0)|^2 \big]^\frac{1}{2}}{\dkl{\w}{x}{a}} \nonumber \\
&\leq \frac{L}{\dkl{\w}{x}{a}} \E \bigg[ \sup_{b \in A} \bigg\{ \frac{| \dkl{v}{X^{x,a}}{b}|}{\dkl{\w}{X^{x,a}}{b}} \dkl{\w}{X^{x,a}}{b} \bigg\}^2 \bigg]^\frac{1}{2} + c_f \nonumber \\
&\leq \frac{L}{\dkl{\w}{x}{a}} \sup_{(y,b) \in \X \times A} \bigg\{ \frac{|\dkl{v}{y}{b}|}{\dkl{\w}{y}{b}} \bigg\} \E \big[ \sup_{b \in A} \w(X^{x,a})(b)^2 \big]^\frac{1}{2} + c_f \nonumber \\
&\leq c_\w L \sup_{(y,b) \in \X \times A} \bigg\{ \frac{|\dkl{v}{y}{b}|}{\dkl{\w}{y}{b}} \bigg\} + c_f.
\end{align}
Combining this, the assumption that $\sup_{(x,a)\in \X \times A} [|\dkl{\w}{x}{a}|^{-1} |\dkl{v}{x}{a}|] < \infty$, and the assumption that $c_\w L < 1$ shows that 
\begin{align}
\sup_{(x,a)\in \X \times A} \frac{|\dkl{v}{x}{a}|}{\dkl{\w}{x}{a}} \leq \frac{c_f}{1- c_\w L}.
\end{align}
This and the assumption that for all $(x,a) \in \X \times A$ it holds that $\E \big[ |\sup_{b \in A} \dkl{\w}{X^{x,a}}{b}|^2 \big]^\frac{1}{2} \leq c_\w \dkl{\w}{x}{a}$ imply that for all $(x,a) \in \X \times A$ it holds that
\begin{align}
\frac{\big( \E\big[ \sup_{b \in A} |\dkl{v}{X^{x,a}}{b}|^2 \big] \big)^\frac{1}{2}}{\dkl{\w}{x}{a}} &\leq \bigg[\sup_{(y,b) \in \X \times A} \frac{|\dkl{v}{y}{b}|}{\dkl{\w}{y}{b}} \bigg] \frac{\E \big[ \sup_{b \in A} \dkl{\w}{X^{x,a}}{b}^2 \big]^\frac{1}{2}}{\dkl{\w}{x}{a}} \leq \frac{c_f c_\w}{1- c_\w L}.
\end{align}
Taking the supremum over $\X \times A$ yields
\begin{align}
\sup_{(x,a) \in \X \times A} \frac{\big( \E\big[ \sup_{b \in A} |\dkl{v}{X^{x,a}}{b}|^2 \big] \big)^\frac{1}{2}}{\dkl{\w}{x}{a}} \leq \frac{c_f c_\w}{1- c_\w L}.
\end{align}
This completes the proof of \cref{sol_sto_fp_eq_estimates}. 
\end{mproof}

\section{Full-history recursive multilevel fixed-point (MLFP) approximations} \label{sec:mlfp_approx}

\subsection{Mathematical description of MLFP approximations}

\begin{setting}\label{mlp_setting}
Let $M \in \N$, let $\Theta = \bigcup_{n \in \N} \Z^n$, let $(\X, \mathcal{X})$ be a nonempty measurable space, let $(\Omega, \mathcal{F}, \mathbb{P})$ be a probability space, let $A$ be a nonempty set, let $\R^A = \{ r \colon A \rightarrow \R \}$, let $\A = \bigotimes_{a \in A} \mathcal{B}(\R)$, let $\left\| \cdot \right\|_\infty \colon \R^A \rightarrow [0,\infty]$ satisfy for all $r \in \R^A$ that $\left\| r \right\|_\infty = \sup_{a \in A} \left| r(a) \right|$, let $\X^A = \{ q \colon A \rightarrow \X \}$, let $\wt{\mathcal{X}} = \bigotimes_{a \in A} \mathcal{X}$, let $f\colon \X \times \R^A \rightarrow \R$ be $(\mathcal{X} \otimes \A)/\mathcal{B}(\R)$-measurable, let $(\mathcal{F}^\theta)_{\theta \in \Theta}$, with $\mathcal{F}^\theta \subseteq \mathcal{F}$, $\theta \in \Theta$, be independent $\sigma$-algebras, let $X^\theta = \big( X^{\theta, x, a}(\omega) \big)_{x \in \X,\; a \in A,\; \omega \in \Omega} \colon \X \times \Omega \rightarrow \X^A$, $\theta \in \Theta$, be i.i.d.\ random fields, such that for all $\theta\in \Theta$ it holds that $X^\theta$ is $(\mathcal{X} \otimes \mathcal{F}^\theta)/\wt{\mathcal{X}}$-measurable, let $V_n^\theta \colon \X \times \Omega \rightarrow \R^A$, $n \in \N_0$, $\theta \in \Theta$, satisfy for all $n \in \N_0$, $\theta \in \Theta$, $x \in \X$, $a \in A$ that 
\begin{align}\label{mlp_approx_v_n}
\dkl{V_n^\theta}{x}{a} = \sum_{ l = 0}^{n-1} \frac{1}{M^{n - l}} \sum_{i = 1}^{M^{n-l}}\Big[ &f\big( X^{(\theta,l,i),x,a}, V^{(\theta,l,i)}_l(X^{(\theta, l ,i),x,a}) \big)\nonumber\\
&- \mathbbm{1}_{\N}(l) f\big( X^{(\theta,l,i),x,a}, V^{(\theta,-l,i)}_{\max\{l - 1, 0\}}(X^{(\theta, l ,i),x,a}) \big) \Big].
\end{align}
\end{setting}

\subsection{Measurability and distributional properties for MLFP approximations}

\begin{lemma}\label{mlp_approx_meas}
Assume \cref{mlp_setting}. It holds for all $n \in \N_0$, $\theta \in \Theta$ that $V_n^\theta$ is $(\mathcal{X} \otimes \sigma( \bigcup_{\eta \in \Theta} \mathcal{F}^{(\theta , \eta)} ))/\A$-measurable and it holds for all $n \in \N_0$, $\theta, \vartheta \in \Theta$ that 
\begin{align}
\X \times \Omega \ni (x,\omega) \mapsto \big[ A \ni a \mapsto f\big( X^{\theta, x, a}(\omega), V^\vartheta_n(X^{\theta , x, a}(\omega), \omega) \big) \in \R \big] \in \R^A
\end{align}
is $(\mathcal{X} \otimes \sigma(\bigcup_{\eta \in \Theta} \mathcal{F}^{(\vartheta, \eta)} \cup \mathcal{F}^\theta))/\A$-measurable.
\end{lemma}

\begin{mproof}{\cref{mlp_approx_meas}}
Note that (\ref{mlp_approx_v_n}) implies that for all $\theta \in \Theta$, $x \in \X$, $a \in A$ it holds that $\dkl{V_0^\theta}{x}{a} = 0$. Hence for all $\theta \in \Theta$ it holds that $V_0^\theta$ is $(\mathcal{X} \otimes \sigma \big( \bigcup_{\eta \in \Theta} \mathcal{F}^{(\theta, \eta)} \big))/\A$-measurable. Fix $n \in \N$ and assume for all $l \in \{0,1,\dots, n-1\}$, $\theta \in \Theta$ that $V_l^\theta$ is $(\mathcal{X} \otimes \sigma \big( \bigcup_{\eta \in \Theta} \mathcal{F}^{(\theta, \eta)} \big))/\A$-measurable. Observe that for all $l \in \{ 0,1,\dots , n-1\}$, $\theta \in \Theta$ it holds that 
\begin{align}\label{mlp_approx_meas_outer}
\X \times \Omega \ni (x,\omega) \mapsto (x, V_l^\theta(x,\omega)) \in \X \times \R^A
\end{align}
is $(\mathcal{X} \otimes \sigma \big( \bigcup_{\eta \in \Theta} \mathcal{F}^{(\theta, \eta)} \big)) / (\mathcal{X} \otimes \A)$-measurable. Since for all $\theta \in \Theta$ it holds that $X^\theta$ is $(\mathcal{X} \otimes \mathcal{F}^\theta)/\wt{\mathcal{X}}$-measurable it follows that for all $\theta \in \Theta$, $a \in A$ it holds that $\X \times \Omega \ni (x, \omega) \mapsto X^{\theta, x, a}(\omega) \in \X$ is $(\mathcal{X} \otimes \mathcal{F}^\theta)/\mathcal{X}$-measurable. This implies that for all $\theta,\vartheta \in \Theta$, $a \in A$ it holds that $\X \times \Omega \ni (x,\omega) \mapsto (X^{\theta,x,a}(\omega),\omega) \in \X \times \Omega$ is $(\mathcal{X} \otimes \sigma \big( \bigcup_{\eta \in \Theta} \mathcal{F}^{(\vartheta, \eta)} \cup \mathcal{F}^\theta \big))/(\mathcal{X} \otimes \sigma\big( \bigcup_{\eta \in \Theta} \mathcal{F}^{(\vartheta, \eta)} \big))$-measurable. This and (\ref{mlp_approx_meas_outer}) ensure that for all $l \in \{ 0,1,\dots , n-1\}$, $\theta, \vartheta \in \Theta$, $a \in A$ it holds that $\X \times  \Omega \ni (x,\omega) \mapsto \big( X^{\theta, x,a}(\omega), V_l^\vartheta(X^{\theta,x,a}(\omega), \omega ) \big) \in \X \times \R^A$ is $(\mathcal{X} \otimes \sigma \big( \bigcup_{\eta \in \Theta} \mathcal{F}^{(\vartheta, \eta)} \cup \mathcal{F}^\theta \big))/(\mathcal{X} \otimes \A)$-measurable. Therefore for all $l \in \{ 0,1,\dots,n-1\}$, $\theta, \vartheta \in \Theta$, $a \in A$ it holds that $\X \times \Omega \ni (x,\omega) \mapsto f \big(  X^{\theta, x,a}(\omega), V_l^\vartheta(X^{\theta,x,a}(\omega), \omega ) \big) \in \R$ is $(\mathcal{X} \otimes \sigma \big( \bigcup_{\eta \in \Theta} \mathcal{F}^{(\vartheta, \eta)} \cup \mathcal{F}^\theta \big))/\mathcal{B}(\R)$-measurable. This implies that for all $l \in \{ 0,1,\dots,n-1 \}$, $\theta, \vartheta \in \Theta$ it holds that
\begin{align}\label{mlp_approx_meas_comp}
\X \times \Omega \ni (x, \omega) \mapsto \big[A \ni a \mapsto f\big(  X^{\theta, x,a}(\omega), V_l^\vartheta(X^{\theta,x,a}(\omega), \omega )\big) \in \R \big] \in \R^A
\end{align}
is $(\mathcal{X} \otimes \sigma \big( \bigcup_{\eta \in \Theta} \mathcal{F}^{(\vartheta, \eta)} \cup \mathcal{F}^\theta \big))/\A$-measurable. Note that for all $l \in \{ 0,1,\dots, n-1\}$, $i \in \N$, $\theta \in \Theta$ it holds that $\sigma \big( \bigcup_{\eta \in \Theta} \mathcal{F}^{(\theta,l,i,\eta)} \cup \mathcal{F}^{(\theta, l, i)} \big) \subseteq \sigma \big( \bigcup_{\eta \in \Theta} \mathcal{F}^{(\theta , \eta)} \big)$ and $\sigma \big( \bigcup_{\eta \in \Theta} \mathcal{F}^{(\theta,-l,i,\eta)} \cup \mathcal{F}^{(\theta, -l, i)} \big) \subseteq \sigma \big( \bigcup_{\eta \in \Theta} \mathcal{F}^{(\theta , \eta)} \big)$. This and (\ref{mlp_approx_meas_comp}) ensure that for all $\theta \in \Theta$ it holds that $V_n^\theta \colon \X \times \Omega \rightarrow \R^A$ is $(\mathcal{X} \otimes \sigma \big( \bigcup_{\eta \in \Theta} \mathcal{F}^{(\theta, \eta)} \big))/\A$-measurable. Induction hence proves that for all $n \in \N_0$, $\theta \in \Theta$ it holds that $V_n^\theta$ is $(\mathcal{X} \otimes \sigma \big( \bigcup_{\eta \in \Theta} \mathcal{F}^{(\theta, \eta)} \big))/\A$-measurable. This and (\ref{mlp_approx_meas_comp}) imply that for all $n \in \N_0$, $\theta, \vartheta \in \Theta$ it holds that $\X \times \Omega \ni (x, \omega) \mapsto \big[  A \ni a \mapsto f\big( X^{\theta, x, a}(\omega), V^\vartheta_n(X^{\theta , x, a}(\omega), \omega) \big) \in \R \big] \in \R^A$ is $(\mathcal{X} \otimes \sigma(\bigcup_{\eta \in \Theta} \mathcal{F}^{(\vartheta, \eta)} \cup \mathcal{F}^\theta))/\A$-measurable. The proof of \cref{mlp_approx_meas} is thus completed.
\end{mproof}

\begin{lemma}\label{mlp_approx_dist}
Assume \cref{mlp_setting}. Then for all $n \in \N_0$ it holds that $V_n^\theta$, $\theta \in \Theta$, are identically distributed random fields.
\end{lemma}

\begin{mproof}{\cref{mlp_approx_dist}}
For all $\theta \in \Theta$, $(x,a) \in \X \times A$ it holds that $\dkl{V_0^\theta}{x}{a} = 0$. Therefore $V_0^\theta, \theta \in \Theta$, are identically distributed random fields. Fix $n \in \N$ and assume that for all $l \in \{0,1,\dots, n-1\}$ it holds that $V_l^\theta$, $\theta \in \Theta$, are identically distributed random fields. \cref{mlp_approx_meas} and \cite[Lemma 2.6]{beck2020nonlinear}
ensure for all $l \in \{ 0,1,\dots , n-1\}$, $i \in \N$, $\theta \in \Theta$ that $\X \times \Omega \ni (x, \omega) \mapsto (V_l^{(\theta, l ,i)}(x,\omega), V_{\max\{ l-1,0 \}}^{(\theta, -l ,i)}(x,\omega)) \in \R^A \times \R^A$ and $\X \times \Omega \ni (x, \omega) \mapsto (V_l^{(0, l ,i)}(x,\omega), V_{\max\{l-1,0\}}^{(0, -l ,i)}(x,\omega) ) \in \R^A \times \R^A$ are identically distributed random fields.
This ensures that for all $l \in \{0,1,\dots , n-1\}$, $i\in \N$, $\theta \in \Theta$ it holds that $\X \times \Omega \ni (x , \omega) \mapsto f(x,V_l^{(\theta, l ,i)}(x,\omega)) - \mathbbm{1}_\N(l)f(x,V_{\max\{l-1,0 \}}^{(\theta, -l ,i)}(x,\omega)) \in \R$ and $\X \times \Omega \ni (x , \omega) \mapsto f(x,V_l^{(0, l ,i)}(x,\omega)) - \mathbbm{1}_\N(l)f(x,V_{\max\{l-1,0 \}}^{(0, -l ,i)}(x,\omega)) \in \R$ are identically distributed random fields. Combining this, the assumption that for all $\theta \in \Theta$ it holds that $X^\theta$ is $(\mathcal{X} \otimes \mathcal{F}^\theta)/\wt{\mathcal{X}}$-measurable, and \cite[Lemma 2.5]{beck2020nonlinear} establishes that for all $l\in \{ 0,1,\dots,n-1\}$, $i \in \N$, $\theta \in \Theta$ it holds that
\begin{align}
\X \times A \times \Omega \ni (x,a,\omega) &\mapsto f\big( X^{(\theta, l ,i),x,a}(\omega), V_l^{(\theta, l ,i)}\big( X^{(\theta, l ,i), x,a}(\omega), \omega\big)\big)\nonumber\\
&\hspace{0.5cm} - \mathbbm{1}_\N(l)f\big( X^{(\theta, l ,i),x,a}(\omega), V_{\max\{l-1,0\}}^{(\theta, -l ,i)}\big( X^{(\theta, l ,i), x,a}(\omega), \omega\big) \big) \in \R
\end{align}
and
\begin{align}
\X \times A \times \Omega \ni (x,a,\omega) &\mapsto f\big( X^{(0, l ,i),x,a}(\omega), V_l^{(0, l ,i)}\big( X^{(0, l ,i), x,a}(\omega), \omega\big)\big) \nonumber \\
&\hspace{0.5cm} - \mathbbm{1}_\N(l)f\big( X^{(0, l ,i),x,a}(\omega), V_{\max\{l-1,0\}}^{(0, -l ,i)}\big( X^{(0, l ,i), x,a}(\omega), \omega\big) \big) \in \R
\end{align}
are identically distributed random fields.
This implies that for all $l \in \{0,1,\dots,n-1\}$, $i \in \N$, $\theta \in \Theta$ it holds that 
\begin{align}
\hspace{-0.3cm}\X \times \Omega \ni (x,\omega) &\mapsto \Big[ A \ni a \mapsto f\big( X^{(\theta, l ,i),x,a}(\omega), V_l^{(\theta, l ,i)}\big( X^{(\theta, l ,i), x,a}(\omega), \omega\big)\big) \nonumber \\
&\hspace{2.2cm} - \mathbbm{1}_\N(l)f\big( X^{(\theta, l ,i),x,a}(\omega), V_{\max\{l-1,0\}}^{(\theta, -l ,i)}\big( X^{(\theta, l ,i), x,a}(\omega), \omega\big) \big) \in \R\Big] \in \R^A
\end{align}
and
\begin{align}
\hspace{-0.3cm}\X \times \Omega \ni (x,\omega) &\mapsto \Big[ A \ni a \mapsto f\big( X^{(0, l ,i),x,a}(\omega), V_l^{(0, l ,i)}\big( X^{(0, l ,i), x,a}(\omega), \omega\big)\big) \nonumber \\
&\hspace{2.2cm} - \mathbbm{1}_\N(l)f\big( X^{(0, l ,i),x,a}(\omega), V_{\max\{l-1,0\}}^{(0, -l ,i)}\big( X^{(0, l ,i), x,a}(\omega), \omega\big) \big) \in \R \Big] \in \R^A
\end{align}
are identically distributed random fields. Let $g_k\colon (\R^A)^k \rightarrow \R^A$, $k \in \N$, satisfy for all $k\in\N$, $r_1, \dots, r_k \in \R^A$ that $g_k(r_1, \dots, r_k) = \sum_{j = 1}^kr_j$. Note that for all $k \in \N$ it holds that $g_k$ is $\A^{\otimes k}/\A$-measurable. Let $U^\theta_{l,i} \colon \X \times \Omega \rightarrow \R^A$, $l \in \{ 0,1,\dots, n-1\}$, $i \in \N$, $\theta \in \Theta$, satisfy for all $l \in \{0,1,\dots,n-1\}$, $i \in \N$, $\theta \in \Theta$, $(x,\omega) \in \X \times \Omega$ that
\begin{align}
U^\theta_{l,i}(x,\omega) &= \big[ A \ni a \mapsto f\big( X^{(\theta, l ,i),x,a}(\omega), V_l^{(\theta, l ,i)}\big( X^{(\theta, l ,i), x,a}(\omega), \omega\big)\big) \nonumber \\
&\hspace{2.4cm} - \mathbbm{1}_\N(l)f\big( X^{(\theta, l ,i),x,a}(\omega), V_{\max\{l-1,0\}}^{(\theta, -l ,i)}\big( X^{(\theta, l ,i), x,a}(\omega), \omega\big) \big) \in \R\big]. 
\end{align}
\cref{mlp_approx_meas} demonstrates that for all $l \in \{0,1,\dots,n-1\}$, $i \in \N$, $\theta \in \Theta$ it holds that $U^\theta_{l,i}$ is $(\mathcal{X} \otimes \sigma ( \bigcup_{\eta \in \Theta} \mathcal{F}^{(\theta, l,i,\eta)} \cup \mathcal{F}^{(\theta, -l,i,\eta)} \cup \mathcal{F}^{(\theta, l ,i)} ) )/\A$-measurable. Note that (\ref{mlp_approx_v_n}) ensures for all $(x,\omega) \in \X \times \Omega$, $\theta \in \Theta$ that $V_n^\theta(x, \omega) = \sum_{l = 0}^{n-1} \frac{1}{M^{n-l}} \sum_{i = 1}^{M^{n-l}} U^\theta_{l,i}(x,\omega)$. Let $\wt{M} = \sum_{j = 1}^nM^j \in \N$, let $Y^\theta_n \colon \X \times \Omega \rightarrow (\R^A)^{\wt{M}}$, $\theta \in \Theta$, satisfy for all $\theta\in \Theta$ that
\begin{align}
Y^\theta_n &= \Big( \frac{1}{M^n} U^\theta_{0,1}, \frac{1}{M^n} U^\theta_{0,2}, \dots, \frac{1}{M^n} U^\theta_{0,M^n} , \frac{1}{M^{n-1}}U^\theta_{1,1}, \frac{1}{M^{n-1}} U^\theta_{1,2},\dots, \frac{1}{M^{n-1}} U^\theta_{1,M^{n-1}} , \dots \nonumber \\
&\hspace{1cm}\dots, \frac{1}{M} U^\theta_{n-1,1}, \frac{1}{M} U^\theta_{n-1,2}, \dots, \frac{1}{M} U^\theta_{n-1,M} \Big).
\end{align}
Observe that for all $(x,\omega) \in \X \times \Omega$, $\theta \in \Theta$ it holds that
\begin{align}
g_{\wt{M}} (Y^\theta_n(x,\omega)) = \sum_{l = 0}^{n-1} \frac{1}{M^{n-l}} \sum_{i = 1}^{M^{n-l}} U^\theta_{l,i}(x,\omega) = V_n^\theta(x, \omega).
\end{align}
The fact that for all $l \in \{ 0,1,\dots, n-1 \}$, $i \in \N$, $\theta \in \Theta$ it holds that $U_{l,i}^\theta$ is $(\mathcal{X} \otimes \sigma ( \bigcup_{\eta \in \Theta} \mathcal{F}^{(\theta, l,i,\eta)} \cup \mathcal{F}^{(\theta, -l,i,\eta)} \cup \mathcal{F}^{(\theta, l ,i)} ) )/\A$-measurable, the assumption that $(\mathcal{F}^\theta)_{\theta \in \Theta}$ are independent $\sigma$-algebras, and \cite[Lemma 2.6]{beck2020nonlinear}
prove for all $\theta \in \Theta$ that $Y^\theta_n$ and $Y^0_n$ are identically distributed random fields.
This implies for all $\theta \in \Theta$ that $V_n^\theta$ and $V_n^0$ are identically distributed random fields. Induction and the fact that $V_0^\theta$, $\theta \in \Theta,$ are identically distributed random fields prove that for all $n \in \N_0$ it holds that $V_n^\theta$, $\theta \in \Theta$, are identically distributed random fields. This completes the proof of \cref{mlp_approx_dist}.
\end{mproof}

\subsection{Integrability properties for MLFP approximations}

\begin{lemma}\label{mlp_approx_2nd_moment} %(2nd moments)(c.f. Lemma 3.4 in stopping)\\
Assume \cref{mlp_setting}, assume $A$ is finite, let $L \in [0,\infty)$, let $\w \colon \X \rightarrow (0, \infty)$ be $\mathcal{X}/ \mathcal{B}((0,\infty))$-measurable, assume that for all $(x,a) \in \X\times A$, $r,s \in \R^A$ it holds that $|f(x,r) - f(x,s)| \leq L \max_{b \in A} |r(b) - s(b)|$, and $\big( \E \big[ |f(X^{0,x,a}, 0 )|^2 + |\w(X^{0,x,a})|^2 \big] \big)^\frac{1}{2} \leq L\w(x)$. Then it holds for all $n \in \N_0$ that $\sup_{x \in \X} \big( \E \big[ \max_{a \in A} \left|\dkl{V_n^0}{x}{a} \right|^2 \big] |\w(x)|^{-2} \big) < \infty $.
\end{lemma}

\begin{mproof}{\cref{mlp_approx_2nd_moment}}
For all $(x,a) \in \X \times A$, $\theta \in \Theta$ it holds that $\dkl{V_0^\theta}{x}{a} = 0$. Hence it holds that $\sup_{x \in \X} (\E \big[ \| V_0^0(x) \|^2_\infty \big] |\w(x)|^{-2} ) < \infty$. Fix $n \in \N$ and assume for all $l \in \{ 0,1,\dots,n-1\}$ that $\sup_{x \in \X} (\E \big[ \| V_l^0(x) \|^2_\infty \big] |\w(x)|^{-2} ) < \infty$. Note that for all $m \in \N$ it holds that
\begin{align} \label{mlp_approx_2nd_moment_eq_1}
\sup_{x \in \X} \frac{\E \big[ \| V_m^0(x) \|^2_\infty \big]}{|\w(x)|^2}  \leq \sup_{x \in \X} \sum_{a \in A} \frac{\E \big[ |\dkl{V_m^0}{x}{a}|^2 \big]}{|\w(x)|^2}  \leq \sum_{a \in A} \sup_{x \in X} \frac{\E \big[ |\dkl{V_m^0}{x}{a}|^2 \big]}{|\w(x)|^2} .
\end{align}
\cref{mlp_approx_dist} and Jensen's inequality prove that for all $m \in \N$, $(x, a) \in \X \times A$, $\theta \in \Theta$ it holds that
\begin{align} \label{mlp_approx_2nd_moment_v_n_jensen}
\E \big[ |\dkl{V_m^0}{x}{a}|^2 \big] &= \E \big[ |\dkl{V_m^\theta}{x}{a}|^2 \big] \nonumber \\
&\leq \sum_{l = 0}^{m-1} \frac{m}{M^{m-l}} \sum_{i = 1}^{M^{m-l}} \E \Big[ \big| f\big( X^{(\theta, l,i), x,a}, V_l^{(\theta,l,i)}(X^{(\theta, l, i),x,a})\big) \nonumber \\
&\hspace{3.5cm}- \mathbbm{1}_\N(l)f\big( X^{(\theta, l, i), x, a}, V_{\max \{ l-1,0\}}^{(\theta, -l, i)}(X^{(\theta, l, i), x, a}) \big) \big|^2 \Big] \nonumber \\
&\leq \sum_{l = 0}^{m-1} \frac{2m}{M^{m-l}} \sum_{i = 1}^{M^{m-l}} \E \big[ |f\big( X^{(\theta, l,i), x,a}, V_l^{(\theta,l,i)}(X^{(\theta, l, i),x,a})\big)|^2 \big] \nonumber \\
&\hspace{3cm}+ \E \big[ |f\big( X^{(\theta, l, i), x, a}, V_{\max \{ l-1,0\}}^{(\theta, -l, i)}(X^{(\theta, l, i), x, a}) \big)|^2 \big].
\end{align}
The triangle inequality, Jensen's inequality, the assumption that for all $x \in \X$, $r, s \in \R^A$ it holds that $|f(x,r) - f(x,s)| \leq L \max_{a \in A} |r(a) - s(a)|$, and the assumption that for all $(x,a) \in \X \times A$ it holds that $\big( \E\big[ |f(X^{0,x,a}, 0)|^2 \big] \big)^\frac{1}{2} \leq L\w(x)$ ensure that for all $m \in \N_0$, $(x,a) \in \X \times A$, $\theta, \vartheta \in \Theta$ it holds that
\begin{align} \label{mlp_approx_2nd_moment_f_decomp} 
\E \big[ |f\big( X^{\theta, x, a}, V_m^\vartheta(X^{\theta, x,a }) \big)|^2 \big] &\leq 2 \E \big[ |f\big( X^{\theta, x, a}, V_m^\vartheta(X^{\theta, x,a }) \big) - f(X^{\theta, x, a}, 0)|^2 \big] + 2 \E \big[ |f(X^{\theta, x, a}, 0)|^2 \big] \nonumber \\
&\leq 2L^2 \E \big[ \max_{b \in A} |\dkl{V_m^\vartheta}{X^{\theta, x, a}}{b}|^2 \big] + 2L^2|\w(x)|^2.
\end{align}
\cref{mlp_approx_dist} ensures for all $m \in \N$, $\vartheta \in \Theta$ that $V_m^\vartheta$ and $V_m^0$ are identically distributed. Combining this and the assumption that $A$ is finite demonstrates that for all $m \in \N$, $\vartheta \in \Theta$ it holds that $\X\times \Omega \ni (x, \omega) \mapsto \| V_m^\vartheta(x,w) \|^2_\infty \in [0,\infty)$ and $\X \times \Omega \ni (x,\omega) \mapsto \| V_m^0(x,\omega) \|^2_\infty \in [0,\infty)$ are identically distributed. This, \cref{mlp_approx_meas}, the assumption that for all $\theta \in \Theta$ it holds that $X^\theta$ is $(\mathcal{X} \otimes \mathcal{F}^\theta)/\wt{\mathcal{X}}$-measurable, the assumption that the $\sigma$-algebras $\mathcal{F}^\theta$, $\theta \in \Theta$, are independent, and
\cite[Lemma 2.2]{HJKNW2018} imply that for all $m \in \N_0$, $(x,a) \in \X \times A$, $\theta, \vartheta \in \Theta$ with $\theta \not \in \bigcup_{\eta \in \Theta} \{(\vartheta, \eta)\}$ it holds that 
\begin{align}
\E \big[ \| V_m^\vartheta(X^{\theta, x,a}) \|^2_\infty \big] &= \int_\X \E \big[ \| V_m^\vartheta(y) \|^2_\infty\big] \big( X^{\theta, x, a} (\mathbb{P}) \big)(dy) \nonumber \\
&= \int_\X \E \big[ \| V_m^0(y) \|^2_\infty \big] \big( X^{0,x,a} (\mathbb{P}) \big)(dy).
\end{align}
This, (\ref{mlp_approx_2nd_moment_f_decomp}), the assumption that for all $(x, a) \in \X \times A$ it holds that $(\E \big[ |\w(X^{0,x,a})|^2 \big])^\frac{1}{2} \leq L\w(x)$, and the assumption that for all $l \in \{0,1,\dots, n-1\}$ it holds that $\sup_{x \in \X}( \E \big[ \| V_l^0(x) \|^2_\infty \big] |\w(x)|^{-2}) < \infty$ prove for all $l \in \{ 0,1, \dots, n-1\}$, $(x,a) \in \X \times A$, $\theta, \vartheta \in \Theta$ with $\theta \not \in \bigcup_{\eta \in \Theta} \{(\vartheta, \eta)\}$ that
\begin{align}  
\frac{ \E \big[ |f\big( X^{\theta, x, a}, V_l^\vartheta(X^{\theta, x,a }) \big)|^2 \big]}{|\w(x)|^2} &\leq \frac{2L^2}{|\w(x)|^2} \int_\X \E \big[ \| V_l^0(y) \|^2_\infty \big] \big( X^{0,x,a} (\mathbb{P}) \big)(dy) + 2L^2 \nonumber \\
&\leq 2L^2 \bigg(\sup_{y \in \X} \frac{\E \big[ \| V_l^0(y) \|^2_\infty \big]}{|\w(y)|^2}\bigg) \bigg( \frac{\E\big[ |\w(X^{0,x,a})|^2 \big]}{|\w(x)|^2} \bigg)  + 2L^2 \nonumber\\
&\leq 2L^4 \bigg(\sup_{y \in \X} \frac{\E \big[ \| V_l^0(y) \|^2_\infty \big]}{|\w(y)|^2}\bigg) + 2L^2.
\end{align}
This ensures that for all $l \in \{ 0,1,\dots,n-1 \}$, $\theta, \vartheta \in \Theta$ with $\theta \not\in \bigcup_{\eta \in \Theta} \{ (\vartheta, \eta) \}$ that
\begin{align}\label{mlp_approx_2nd_moment_esti_f}
\sup_{(x,a) \in \X \times A} \frac{ \E \big[ |f\big( X^{\theta, x, a}, V_l^\vartheta(X^{\theta, x,a }) \big)|^2 \big]}{|\w(x)|^2} &\leq 2L^4 \bigg(\sup_{y \in \X} \frac{\E \big[ \| V_l^0(y) \|^2_\infty \big]}{|\w(y)|^2}\bigg) + 2L^2 < \infty.
\end{align}
This and (\ref{mlp_approx_2nd_moment_v_n_jensen}) demonstrate for all $a \in A$ that $\sup_{x \in \X} \big( \E \big[ |\dkl{V_n^0}{x}{a}|^2 \big] |\w(x)|^{-2} \big) < \infty$. This, the assumption that $A$ is finite, and (\ref{mlp_approx_2nd_moment_eq_1}) show that $\sup_{x \in \X} \big(\E \big[ \| V_n^0(x) \|^2_\infty \big] |\w(x)|^{-2} \big) < \infty$. Induction hence proves that for all $n \in \N_0$ it holds that $\sup_{x \in \X} \big( \E \big[ \max_{a \in A} | \dkl{V_n^0}{x}{a} |^2 \big] |\w(x)|^{-2} \big) < \infty$. The proof of \cref{mlp_approx_2nd_moment} is thus completed.
\end{mproof}

\begin{lemma} \label{mlp_approx_ew} %(Mean of MLP approximations) (c.f. Lemma 3.5 in stopping)\\
Assume \cref{mlp_setting}, assume $A$ is finite, let $L \in [0,\infty)$, let $\w \colon \X \rightarrow (0, \infty)$ be $\mathcal{X}/\mathcal{B}((0,\infty))$-measurable, assume that for all $(x,a) \in \X \times A$, $r,s \in \R^A$ it holds that $|f(x,r) - f(x,s)| \leq L \max_{b \in A} |r(b) - s(b)|,$ $\big( \E \big[ |f(X^{0,x,a}, 0)|^2 + |\w(X^{0,x,a})|^2  \big] \big)^\frac{1}{2} \leq L \w(x)$. Then it holds for all $n \in \N$, $(x,a) \in \X \times A$, $\theta \in \Theta$ that $\E \big[ |f(X^{0,x,a}, V_{n-1}^0(X^{0,x,a}))| + |\dkl{V_n^\theta}{x}{a}| \big] < \infty$ and $\E \big[ \dkl{V_n^\theta}{x}{a} \big] = \E \big[ f(X^{0,x,a}, V_{n-1}^0(X^{0,x,a})) \big]$.
\end{lemma}

\begin{mproof}{\cref{mlp_approx_ew}}
The assumption that for all $x \in \X$, $r,s \in \R^A$ it holds that $|f(x,r) - f(x,s)| \leq L \max_{a \in A} |r(a) - s(a)|$, \cref{mlp_approx_dist}, \cref{mlp_approx_2nd_moment}, the fact that $0 \not \in \bigcup_{\eta \in \Theta} \{(0,\eta)\}$, the assumption that $\big( \E\ \big[ |f(X^{0,x,a}, 0)|^2 \big] \big)^\frac{1}{2} \leq L\w(x)$, and \cite[Lemma 2.2]{HJKNW2018} prove that for all $n \in \N$, $(x,a) \in \X \times A$, $\theta \in \Theta$ it holds that $\E \big[ |f( X^{0,x,a}, V_{n-1}^0(X^{0,x,a}) )|^2 + |\dkl{V_n^\theta}{x}{a}|^2 \big] < \infty$. This ensures for all $n \in \N$, $(x,a) \in \X \times A$, $\theta \in \Theta$ that $\E \big[ |f( X^{0,x,a}, V_{n-1}^0(X^{0,x,a}) )| + |\dkl{V_n^\theta}{x}{a}| \big] < \infty$. This, \cref{mlp_approx_meas}, the assumption that $X^\theta$, $\theta \in \Theta$, are i.i.d.\ random fields, \cref{mlp_approx_dist}, and \cite[Lemma 2.2]{HJKNW2018} demonstrate that for all $n \in \N_0$, $(x,a) \in \X \times A$, $\theta, \vartheta \in \Theta$ with $\theta \not\in \bigcup_{\eta \in \Theta} \{(\vartheta, \eta)\}$ it holds that
\begin{align}
\E \big[ f\big( X^{\theta, x,a}, V_n^\vartheta(X^{\theta ,x,a}) \big) \big] &= \E \big[ f\big( X^{0,x,a}, V_n^0(X^{0,x,a}) \big) \big].
\end{align}
This establishes for all $n \in \N$, $(x,a) \in \X \times A$, $\theta \in \Theta$ that 
\begin{align}
\E \big[ \dkl{V_n^\theta}{x}{a} \big] &= \sum_{l = 0}^{n-1} \frac{1}{M^{n-l}} \sum_{i = 1}^{M^{n-l}} \E \big[ f\big( X^{(\theta, l, i),x,a}, V_l^{(\theta, l, i)}(X^{(\theta,l,i),x,a}) \big) \big] \nonumber \\
&\hspace{3.2cm}- \mathbbm{1}_{\N}(l) \E \big[ f\big( X^{(\theta, l, i),x,a} V_{\max \{l-1,0 \}}^{(\theta, -l ,i)}(X^{(\theta, l ,i)}) \big) \big] \nonumber \\
&= \sum_{l = 0}^{n-1} \E \big[ f\big( X^{0,x,a}, V_l^0(X^{0,x,a}) \big) \big] - \mathbbm{1}_{\N}(l) \E \big[ f \big( X^{0,x,a}, V_{\max\{ l-1,0 \}}^0(X^{0,x,a}) \big) \big] \nonumber \\
&= \E \big[ f \big( X^{0,x,a}, V_{n-1}^0(X^{0,x,a}) \big) \big].
\end{align}
The proof of \cref{mlp_approx_ew} is thus completed. 
\end{mproof}

\subsection{Recursive error bounds for MLFP approximations}

\begin{lemma} \label{mlp_approx_bias}%(Bias)(c.f. Lemma 3.6 in stopping)\\
Assume \cref{mlp_setting}, assume $A$ is finite, let $c, L \in [0,\infty)$, let $\w \colon \X \rightarrow (0, \infty)$ be $\mathcal{X}/\mathcal{B}((0,\infty))$-meaurable, let $v \colon \X \rightarrow \R^A$ be $\mathcal{X}/\A$-measurable, assume that for all $(x,a) \in \X \times A$, $r,s \in \R^A$ it holds that $|f(x,r) - f(x,s)| \leq L \max_{b \in A} |r(b) - s(b)|$, $\big( \E \big[ |f(X^{0,x,a}, 0)|^2 + \|v(X^{0,x,a}) \|_\infty^2 + \left| \w(X^{0,x,a}) \right|^2 \big] \big)^\frac{1}{2} \leq c \w(x)$, and $\dkl{v}{x}{a} = \E [ f(X^{0,x,a}, v(X^{0,x,a})) ]$. Then it holds for all $n \in \N$, $(x,a) \in \X \times A$, $\theta \in \Theta$ that
\begin{align}\label{mlp_approx_bias_1}
\big| \dkl{v}{x}{a} - \E\big[ \dkl{V_n^\theta}{x}{a} \big] \big| \leq L \left( \E \Big[ \big\| v(X^{0,x,a}) - V_{n-1}^0(X^{0,x,a}) \big\|^2_\infty \Big] \right)^{\frac{1}{2}}.
\end{align}
\end{lemma}

\begin{mproof}{\cref{mlp_approx_bias}}
The assumption that for all $(x,a) \in \X \times A$ it holds that $\dkl{v}{x}{a} = \E \big[ f(X^{0,x,a}, v(X^{0,x,a})) \big]$, \cref{mlp_approx_ew}, the assumption that for all $x \in \X$, $r,s \in \R^A$ it holds that $\left| f(x,r) - f(x,s) \right| \leq L \max_{a \in A} \left| r(a) - s(a) \right|$, \cref{mlp_approx_2nd_moment}, and the Cauchy-Schwarz inequality imply that for all $n \in \N$, $(x,a) \in \X \times A$, $\theta \in \Theta$ it holds that
\begin{align} \label{mlp_approx_bias_upto_disint}
\left| v(x)(a) - \E \big[ V_n^\theta(x)(a)\big] \right| &\leq \E \big[ \big| f\big( X^{0,x,a}, v(X^{0,x,a}) \big) - f\big( X^{0,x,a}, V_{n-1}^0(X^{0,x,a}) \big) \big| \big] \nonumber \\
&\leq L \E \Big[ \big\| v(X^{0,x,a})-V_{n-1}^0(X^{0,x,a}) \big\|_\infty \Big] \nonumber \\
&\leq L \Big( \E \Big[ \big\| v(X^{0,x,a}) - V_{n-1}^0(X^{0,x,a}) \big\|^2_\infty \Big] \Big)^\frac{1}{2}. 
\end{align}
The proof of \cref{mlp_approx_bias} is thus completed.
\end{mproof}

\begin{lemma} \label{mlp_approx_var} %(Variance)(c.f. Lemma 3.7 in stopping)\\ 
Assume \cref{mlp_setting}, assume $A$ is finite, let $c,L \in [0, \infty)$, let $\w \colon \X \rightarrow (0,\infty)$ be $\mathcal{X}/\mathcal{B}((0,\infty))$-measurable, let $v \colon \X \rightarrow \R^A$ be $\mathcal{X}/\A$-measurable, assume that for all $(x,a)\in \X \times A$, $r,s \in \R^A$ it holds that $|f(x,r) - f(x,s)| \leq L \max_{b \in A} |r(b) - s(b)|$, $\big( \E \big[ |f(X^{0,x,a}, 0)|^2 + \|v(X^{0,x,a}) \|_\infty^2 + \left|\w(X^{0,x,a})\right|^2\big] \big)^\frac{1}{2} \leq c\w(x)$, and $\dkl{v}{x}{a} = \E \big[ f(X^{0,x,a}, v(X^{0,x,a})) \big]$. Then it holds for all $n \in \N$, $(x,a) \in \X \times A$, $\theta \in \Theta$ that 
\begin{align}
&\hspace{-0.3cm} \left( \var \left[ \dkl{V_n^\theta}{x}{a} \right] \right)^\frac{1}{2} \nonumber \\
&\hspace{-0.3cm} \leq \frac{1}{\sqrt{M^n}} \Big( \E \big[ |f( X^{0,x,a},0 )|^2 \big] \Big)^\frac{1}{2} + \mathbbm{1}_{[2, \infty)}(n) \frac{L}{\sqrt{M}} \big( \E \big[ \| V_{n-1}^0(X^{0,x,a}) - v(X^{0,x,a}) \|^2_\infty \big] \big)^\frac{1}{2} \nonumber \\
&\hspace{-0.1cm}+ \mathbbm{1}_{[2, \infty)}(n) \frac{L\sqrt{M}}{\sqrt{M^n}} \big( \E \big[ \| v(X^{0,x,a}) \|^2_\infty \big] \big)^\frac{1}{2} + \sum_{l = 1}^{n - 2} \frac{L(1 + \sqrt{M})}{\sqrt{M^{n-l}}} \big( \E \big[ \| V_l^0(X^{0,x,a}) - v(X^{0,x,a}) \|^2_\infty \big)^\frac{1}{2}.
\end{align}
\end{lemma}

\begin{mproof}{\cref{mlp_approx_var}}
\cref{mlp_approx_2nd_moment} ensures for all $n \in \N$, $(x,a) \in \X \times A$, $\theta \in \Theta$ that $\var\big[ \dkl{V_n^\theta}{x}{a} \big] < \infty$. \cref{mlp_approx_meas}
yields that for all $l,i \in \N_0$, $(x,a) \in \X \times A$, $\theta \in \Theta$ it holds that
\begin{align}
\Omega \ni \omega &\mapsto f \big( X^{(\theta, l, i), x, a}(\omega), V_l^{(\theta, l, i)}(X^{(\theta, l, i), x, a}(\omega), \omega) \big) \nonumber \\
&\hspace{3cm}-\mathbbm{1}_\N(l) f\big( X^{(\theta, l, i),x,a}(\omega), V_{\max\{l-1,0\}}^{(\theta, -l, i)}(X^{(\theta, l, i),x ,a}(\omega), \omega) \big) \in \R
\end{align}
is $\sigma \big( \bigcup_{\eta \in \Theta} \mathcal{F}^{(\theta, l, i, \eta)} \cup \mathcal{F}^{(\theta, -l, i, \eta)} \cup \mathcal{F}^{(\theta, l, i)} \big) / \mathcal{B}(\R)$-measurable.
This demonstrates that for all $n \in \N_{0}$, $(x,a) \in \X \times A$ it holds that
\begin{align} \label{mlp_approx_var_eq_1_sum_of_var}
\var\left[ \dkl{V_n^\theta}{x}{a} \right] &= \sum_{l = 0}^{n-1} \frac{1}{(M^{n-l})^2} \sum_{i = 1}^{M^{n-l}} \var \Big[ f\big( X^{(\theta, l,i),x,a}, V_l^{(\theta, l,i)}(X^{(\theta, l, i),x,a}) \big)\nonumber\\
&\hspace{4.5cm} - \mathbbm{1}_\N(l) f\big( X^{(\theta, l ,i),x,a}, V_{\max\{ l-1,0 \}}^{(\theta, -l ,i)}(X^{(\theta, l ,i), x ,a}) \big) \Big].
\end{align}
\cref{mlp_approx_meas}, \cref{mlp_approx_dist}, and \cite[Lemma 2.6]{beck2020nonlinear} ensure that for all $l,i \in \N_{0}$, $\theta \in \Theta$ it holds that 
\begin{align}
\X \times \Omega \ni (x, \omega) &\mapsto \big( V_l^{(\theta, l, i)}(x,\omega), V_{\max\{ l-1,0 \}}^{(\theta, -l ,i)}(x,\omega) \big) \in \R^A \times \R^A \quad \text{and} \nonumber \\
\X \times \Omega \ni (x, \omega) &\mapsto \big( V_l^0(x,\omega), V^1_{\max\{ l-1, 0 \}}(x, \omega) \big) \in \R^A \times \R^A.
\end{align}
are identically distributed random fields.
This and the fact that $f$ is $(\mathcal{X} \otimes \A)/\mathcal{B}(\R)$-measurable show that for all $l,i \in \N_{0}$, $\theta \in \Theta$ it holds that
\begin{align}
\X \times \Omega \ni (x, \omega) &\mapsto f(x, V_l^{(\theta, l,i)}(x,\omega)) - \mathbbm{1}_\N(l) f( x, V_{\max\{ l-1, 0 \}}^{(\theta, -l, i)}(x,\omega) ) \in \R \quad \text{and} \nonumber \\
\X \times \Omega \ni (x, \omega) &\mapsto f( x, V_l^0(x, \omega ) ) - \mathbbm{1}_\N(l) f(x , V_{\max\{ l-1,0 \}}^1(x,\omega)) \in \R
\end{align}
are identically distributed random fields. \cite[Lemma 2.5]{beck2020nonlinear} proves that for all $l,i \in \N_0$, $\theta \in \Theta$ it holds that 
\begin{align}
\X \times A \times \Omega \ni (x,a,\omega) &\mapsto f\big( X^{(\theta, l, i), x, a}(\omega), V_l^{(\theta, l, i)}( X^{(\theta, l, i),x,a}(\omega) , \omega) \big) \nonumber \\
&\hspace{1cm}- \mathbbm{1}_\N(l) f\big( X^{(\theta, l, i), x, a}(\omega), V_{\max \{ l-1, 0 \}}^{(\theta, -l,i)} (X^{(\theta, l, i), x, a}(\omega), \omega) \big) \in \R \quad \text{and} \nonumber \\
\X \times A \times \Omega \ni (x,a,\omega) &\mapsto f\big( X^{0, x, a}(\omega), V_l^0( X^{0,x,a}(\omega) , \omega) \big) \nonumber \\
&\hspace{1cm}- \mathbbm{1}_\N(l) f\big( X^{0, x, a}(\omega), V_{\max \{ l-1, 0 \}}^1 (X^{0, x, a}(\omega), \omega) \big) \in \R
\end{align} 
are identically distributed random fields. This implies for all $l,i \in \N_{0}$, $(x,a) \in \X \times A$, $\theta \in \Theta$ that
\begin{align} \label{mlp_approx_var_eq_2} 
&\var \Big[ f\big( X^{(\theta, l,i),x,a}, V_l^{(\theta, l,i)}(X^{(\theta, l, i),x,a}) \big) - \mathbbm{1}_\N(l) f\big( X^{(\theta, l ,i),x,a}, V_{\max\{ l-1,0 \}}^{(\theta, -l ,i)}(X^{(\theta, l ,i), x ,a}) \big) \Big]\nonumber\\
&\hspace{0.5cm} = \var \Big[ f\big( X^{0,x,a}, V_l^0(X^{0,x,a}) \big) - \mathbbm{1}_\N(l) f\big( X^{0,x,a}, V_{\max\{ l-1,0 \}}^1(X^{0, x ,a}) \big) \Big].
\end{align}
This and (\ref{mlp_approx_var_eq_1_sum_of_var})
ensure that for all $n \in \N$, $(x, a) \in \X \times A$, $\theta \in \Theta$ it holds that
\begin{align}
\hspace{-0.3cm} \var[V_n^\theta(x)(a)] &= \sum_{l = 0}^{n-1} \frac{1}{M^{n-l}} \var \Big[ f\big( X^{0,x,a}, V_l^0(X^{0,x,a}) \big) - \mathbbm{1}_\N(l) f\big( X^{0,x,a}, V_{\max\{ l-1,0 \}}^1(X^{0, x ,a}) \big) \Big] \nonumber \\
&\leq \sum_{l = 0}^{n-1} \frac{1}{M^{n-l}} \E \Big[ \big| f\big( X^{0,x,a}, V_l^0(X^{0,x,a}) \big) - \mathbbm{1}_\N(l) f\big( X^{0,x,a}, V_{\max\{ l-1,0 \}}^1(X^{0, x ,a}) \big) \big|^2 \Big].
\end{align}
Combining this, the fact that for all $(x, a) \in \X \times A$ it holds that $ \dkl{V_0^0}{x}{a} = 0$, and the fact that for all $n \in \N$, $r_1, r_2, \dots, r_n \in [0,\infty)$ it holds that $\sqrt{r_1 + r_2 + \dots + r_n} \leq \sqrt{r_1} + \sqrt{r_2} + \dots + \sqrt{r_n}$ prove that for all $n \in \N$, $(x,a) \in \X \times A$, $\theta \in \Theta$ it holds that
\begin{align}
\big( \var\big[&\dkl{V_n^\theta}{x}{a} \big] \big)^\frac{1}{2} \nonumber \\
&\leq \sum_{l = 0}^{n-1} \frac{1}{\sqrt{M^{n-l}}} \Big( \E \Big[ \big| f\big( X^{0,x,a}, V_l^0(X^{0,x,a}) \big) - \mathbbm{1}_\N(l) f\big( X^{0,x,a}, V_{\max\{ l-1, 0 \}}^1(X^{0,x,a}) \big) \big|^2 \Big] \Big)^\frac{1}{2} \nonumber \\
&= \sum_{l = 1}^{n-1} \frac{1}{\sqrt{M^{n-l}}} \Big( \E \Big[ \big| f\big( X^{0,x,a}, V_l^0(X^{0,x,a}) \big) - \mathbbm{1}_\N(l) f\big( X^{0,x,a}, V_{l-1}^1(X^{0,x,a}) \big) \big|^2 \Big] \Big)^\frac{1}{2} \nonumber \\
&\hspace{0.5cm}+ \frac{1}{\sqrt{M^n}} \Big( \E \big[ |f( X^{0,x,a},0 )|^2 \big] \Big)^\frac{1}{2}.
\end{align} 
The assumption that for all $x \in \X$, $r, s \in \R^A$ it holds that $| f(x,r) - f(x,s) | \leq L \max_{a \in A} | r(a) - s(a) |$ and the triangle inequality imply that for all $n \in \N$, $(x,a) \in \X \times A$, $\theta \in \Theta$ it holds that
\begin{align} \label{mlp_approx_var_eq3}
\big( \var\big[& \dkl{V_n^\theta}{x}{a} \big] \big)^\frac{1}{2} \nonumber \\
&\leq \frac{1}{\sqrt{M^n}} \Big( \E \big[ |f( X^{0,x,a},0 )|^2 \big] \Big)^\frac{1}{2} + \sum_{l = 1}^{n-1} \frac{L}{\sqrt{M^{n-l}}} \Big( \E \Big[ \big\| V_l^0(X^{0,x,a}) - V_{l-1}^1(X^{0,x,a}) \big\|^2_\infty \Big] \Big)^\frac{1}{2} \nonumber \\
&\leq \sum_{l = 1}^{n-1} \frac{L}{\sqrt{M^{n-l}}} \bigg[ \big( \E \big[ \| V_l^0(X^{0,x,a}) - v(X^{0,x,a}) \|_\infty^2 \big] \big)^\frac{1}{2} + \big( \E \big[ \| v(X^{0,x,a}) - V_{l-1}^1(X^{0,x,a}) \|_\infty^2 \big] \big)^\frac{1}{2} \bigg] \nonumber \\
&\hspace{0.5cm}+ \frac{1}{\sqrt{M^n}} \Big( \E \big[ |f( X^{0,x,a},0 )|^2 \big] \Big)^\frac{1}{2}.
\end{align} 
The fact that for all $l \in \N_0$, $(x,a) \in \X \times A$ it holds that $\Omega \ni \omega \mapsto V_l^0(X^{0,x,a}(\omega),\omega) \in \R^A$ and $\Omega \ni \omega \mapsto V_l^1(X^{0,x,a}(\omega), \omega) \in \R^A$ are identically distributed ensures for all $l \in \N_0$, $(x,a) \in \X \times A$ that $\Omega \ni \omega \mapsto \| v(X^{0,x,a}(\omega)) - V_l^0(X^{0,x,a}(\omega), \omega) \|^2_\infty \in [0, \infty)$ and $\Omega \ni \omega \mapsto \| v(X^{0,x,a}(\omega)) - V_l^1(X^{0,x,a}(\omega), \omega) \|^2_\infty \in [0, \infty)$ are identically distributed. This and (\ref{mlp_approx_var_eq3}) establish that for all $n \in \N$, $(x,a) \in \X \times A$, $\theta \in \Theta$ it holds that
\begin{align}
&\hspace{-0.3cm}\big( \var\big[ \dkl{V_n^\theta}{x}{a} \big] \big)^\frac{1}{2} \nonumber \\
&\hspace{-0.3cm}\leq \sum_{l = 1}^{n-1} \frac{L}{\sqrt{M^{n-l}}} \bigg[ \big( \E \big[ \| V_l^0(X^{0,x,a}) - v(X^{0,x,a}) \|_\infty^2 \big] \big)^\frac{1}{2} + \big( \E \big[ \| v(X^{0,x,a}) - V_{l-1}^0(X^{0,x,a}) \|_\infty^2 \big] \big)^\frac{1}{2} \bigg] \nonumber \\
&\hspace{-0.1cm}+ \frac{1}{\sqrt{M^n}} \Big( \E \big[ |f( X^{0,x,a},0 )|^2 \big] \Big)^\frac{1}{2} \nonumber \\
&\hspace{-0.3cm}= \frac{1}{\sqrt{M^n}} \Big( \E \big[ |f( X^{0,x,a},0 )|^2 \big] \Big)^\frac{1}{2} + \sum_{l = 1}^{n-1} \frac{L}{\sqrt{M^{n-l}}}\big( \E \big[ \| V_l^0(X^{0,x,a}) - v(X^{0,x,a}) \|_\infty^2 \big] \big)^\frac{1}{2} \nonumber \\
&\hspace{-0.1cm}+ \sum_{l = 0}^{n-2} \frac{L\sqrt{M}}{\sqrt{M^{n-l}}} \big( \E \big[ \| V_l^0(X^{0,x,a}) - v(X^{0,x,a}) \|_\infty^2 \big] \big)^\frac{1}{2} \nonumber \\
&\hspace{-0.3cm}= \frac{1}{\sqrt{M^n}} \Big( \E \big[ |f( X^{0,x,a},0 )|^2 \big] \Big)^\frac{1}{2} + \mathbbm{1}_{[2, \infty)}(n) \frac{L}{\sqrt{M}} \big( \E \big[ \| V_{n-1}^0(X^{0,x,a}) - v(X^{0,x,a}) \|^2_\infty \big] \big)^\frac{1}{2} \nonumber \\
&\hspace{-0.1cm}+ \mathbbm{1}_{[2, \infty)}(n) \frac{L\sqrt{M}}{\sqrt{M^n}} \big( \E \big[ \| v(X^{0,x,a}) \|^2_\infty \big] \big)^\frac{1}{2} + \sum_{l = 1}^{n - 2} \frac{L(1 + \sqrt{M})}{\sqrt{M^{n-l}}} \big( \E \big[ \| V_l^0(X^{0,x,a}) - v(X^{0,x,a}) \|^2_\infty \big)^\frac{1}{2}.
\end{align}
The proof of \cref{mlp_approx_var} is thus completed. 
\end{mproof}

\subsection{Non-recursive error bounds for MLFP approximations}

\begin{propo} \label{mlp_approx_full_error} %(Full Error)(c.f. Proposition 3.8 in stopping)\\
Assume \cref{mlp_setting}, assume $A$ is finite, let $c_f, c_v, c_\w, L \in [0, \infty)$, let $\w \colon \X \rightarrow (0, \infty)$ be $\mathcal{X} / \mathcal{B}((0,\infty))$-measurable, let $v \colon \X \rightarrow \R^A$ be $\mathcal{X}/\A$-measurable, let $c = \frac{3}{2} \max \big\{ \frac{c_v}{c_\w}, c_vL + c_f, \frac{|A|c_f}{c_\w|A|L + 1} \big\}$, assume that for all $(x,a) \in \X \times A$, $r,s \in \R^A$ it holds that $|f(x,r) - f(x,s)| \leq L \max_{b \in A} |r(b) - s(b)|$,\hspace{0.5cm} $\big( \E \big[ |f(X^{0,x,a}, 0)|^2 \big] \big)^\frac{1}{2} \leq c_f \w(x)$, \hspace{0.5cm} $\big( \E \big[ \|v(X^{0,x,a}) \|_\infty^2 \big] \big)^\frac{1}{2} \leq c_v \w(x)$, $\big( \E \big[ \big|\w (X^{0,x,a})\big|^2 \big] \big)^\frac{1}{2} \leq c_\w \w(x)$, and $\dkl{v}{x}{a} = \E \big[ f(X^{0,x,a}, v(X^{0,x,a})) \big]$. Then it holds for all $n \in \N$, $x\in \X$ that 
\begin{align}
\hspace{-0.5cm}\bigg( \frac{\E \big[ \| v(x) - V_n^\theta(x) \|_\infty^2 \big]}{|\w(x)|^2} \bigg)^{\nicefrac{1}{2}} \hspace{-0.2cm} \leq c \Bigg( \tfrac{c_\w L(1 + |A|M^{-\frac{1}{2}}) + M^{-\frac{1}{2}} + \sqrt{\big( c_\w L(1 + |A|M^{-\frac{1}{2}}) + M^{-\frac{1}{2}}\big)^2 + 4c_\w LM^{-\frac{1}{2}}(|A| - 1)}  }{2}  \Bigg)^n \hspace{-0.2cm}.  
\end{align}
\end{propo}

\begin{mproof}{\cref{mlp_approx_full_error}}
The triangle inequality, \cref{mlp_approx_2nd_moment}, and the assumption that $A$ is finite ensure that for all $n \in \N$, $x\in \X$, $\theta \in \Theta$ it holds that
\begin{align} \label{mlp_approx_full_error_base_eq}
\Big( \E \big[ \| v(x) - V_n^\theta(x) \|_\infty^2 \big] \Big)^\frac{1}{2} &= \Big( \E \Big[ \big\| v(x) - \E\big[ V_n^\theta(x) \big] + \E \big[ V_n^\theta(x) \big] - V_n^\theta(x)\big\|_\infty^2 \Big] \Big)^\frac{1}{2}\nonumber\\
&\leq \Big( \E \Big[ \big\| v(x) - \E \big[ V_n^\theta(x) \big] \big\|_\infty^2 \Big] \Big)^\frac{1}{2} + \Big( \E \Big[ \big\| V_n^\theta(x) - \E \big[ V_n^\theta(x) \big] \big\|^2_\infty \Big] \Big)^\frac{1}{2}\nonumber\\
&\leq \big\| v(x) - \E \big[ V_n^\theta(x) \big] \big\|_\infty + \sum_{a \in A} \Big( \E \Big[ \big| V_n^\theta(x)(a) - \E \big[ V_n^\theta(x)(a) \big] \big|^2 \Big] \Big)^\frac{1}{2}\nonumber\\
&= \big\| v(x) - \E \big[ V_n^\theta(x) \big] \big\|_\infty + \sum_{a \in A} \big( \var[V_n^\theta(x)(a)] \big)^\frac{1}{2}.
\end{align}
The assumption that $X^0 \colon \X \times \Omega \rightarrow \X^A$ is $(\mathcal{X} \otimes \mathcal{F}^0)/\wt{\mathcal{X}}$-measurable, the assumption that the $\sigma$-algebras $\mathcal{F}^\theta$, $\theta \in \Theta$, are independent, the assumption that $A$ is finite, \cref{mlp_approx_meas}, and the disintegration theorem \cite[Lemma 2.2]{HJKNW2018} prove that for all $l \in \N$, $(x,a) \in \X \times A$ it holds that
\begin{align} \label{mlp_approx_full_error_disint}
\E \big[ \| V_l^0(X^{0,x,a}) - v(X^{0,x,a}) \|^2_\infty \big] &= \int_\X \E \big[ \| V_l^0(y) - v(y)\|^2_\infty \big] \big( X^{0,x,a}(\mathbb{P}) \big)(dy) \nonumber \\
&\leq \bigg(\sup_{y \in \X} \frac{\E \big[ \| V_l^0(y) - v(y) \|^2_\infty \big]}{|\w(y)|^2} \bigg) c_\w^2 |\w(x)|^2.
\end{align}
This, \cref{mlp_approx_bias}, \cref{mlp_approx_var} and (\ref{mlp_approx_full_error_base_eq}) imply that for all $n \in \N$, $x \in \X$ it holds that
\begin{align} \label{mlp_approx_full_error_eq_2}
\hspace{-1cm}\big(\E& \big[ \| v(x) - V_n^0(x)\|^2_\infty \big]\big)^\frac{1}{2} \nonumber \\
&\leq c_\w \w(x)L \bigg( \sup_{y \in \X} \frac{\E \big[ \| v(y) - V_{n-1}^0(y) \|^2_\infty \big]}{|\w(y)|^2}  \bigg)^\frac{1}{2} + \sum_{a \in A} \bigg[ \frac{1}{\sqrt{M^n}} \big( \E \big[ |f(X^{0,x,a}, 0)|^2 \big] \big)^\frac{1}{2} \nonumber \\
&\hspace{0.5cm}+ \mathbbm{1}_{[2,\infty)(n)} \frac{L\sqrt{M}}{\sqrt{M^n}} \big( \E \big[ \| v(X^{0,x,a}) \|^2_\infty \big] \big)^\frac{1}{2} + \mathbbm{1}_{[2,\infty)}(n) \frac{c_\w \w(x) L}{\sqrt{M}} \bigg( \sup_{y \in \X} \frac{\E \big[ \| v(y) - V_{n-1}^0(y) \|^2_\infty \big]}{|\w(y)|^2}  \bigg)^\frac{1}{2} \nonumber\\
&\hspace{0.5cm}+ \sum_{l = 1}^{n-2}\frac{c_\w \w(x) L(1 + \sqrt{M})}{\sqrt{M^{n-l}}} \bigg( \sup_{y \in X} \frac{\E \big[ \| v(y) - V_l^0(y)\|^2_\infty \big]}{|\w(y)|^2}  \bigg)^\frac{1}{2}\bigg] \nonumber \\
&\leq c_\w \w(x) L \bigg(\sup_{y \in \X} \frac{\E \big[ \| v(y) - V_{n-1}^0(y) \|^2_\infty \big]}{|\w(y)|^2} \bigg)^\frac{1}{2} + \frac{|A|\w(x)}{\sqrt{M^n}} \big(c_f + \mathbbm{1}_{[2,\infty)}(n) c_v L \sqrt{M}\big) \nonumber \\
&\hspace{0.5cm}+ \mathbbm{1}_{[2,\infty)}(n) \frac{c_\w |A|L \w(x)}{\sqrt{M}} \bigg(\sup_{y \in \X} \frac{\E \big[ \| v(y) - V_{n-1}^0(y) \|^2_\infty \big]}{|\w(y)|^2} \bigg)^\frac{1}{2}\nonumber\\
&\hspace{0.5cm}+ c_\w \w(x) |A|L(1 + \sqrt{M}) \sum_{l = 1}^{n-2} \frac{1}{\sqrt{M^{n-l}}} \bigg( \sup_{y \in \X} \frac{\E \big[ \| v(y) - V_l^0(y) \|^2_\infty \big]}{|\w(y)|^2} \bigg)^\frac{1}{2}.
\end{align}
For all $n \in \N_0$ let
\begin{align} \label{mlp_approx_full_error_def_fn_an}
F_n = \bigg( \sup_{x \in \X} \frac{\E \big[ \| v(x) - V_n^0(x) \|^2_\infty \big]}{|\w(x)|^2}  \bigg)^\frac{1}{2} \quad \text{and} \quad a_n &= M^\frac{n}{2} F_n.
\end{align}
\cref{mlp_approx_2nd_moment} ensures that for all $n \in \N_0$ it holds that $a_n, F_n \in [0, \infty)$. Note that (\ref{mlp_approx_full_error_eq_2}) shows for all $n \in \N \cap [2, \infty)$ that
\begin{align}
F_n &\leq |A| M^{- \frac{n}{2}} \big( c_f + c_vLM^\frac{1}{2}\big) + c_\w L(1 + |A|M^{- \frac{1}{2}}) \bigg( \sup_{y \in \X} \frac{\E \big[ \| v(y) - V_{n-1}^0(y) \|^2_\infty\big]}{|\w(y)|^2} \bigg)^\frac{1}{2} \nonumber \\
&\hspace{0.5cm}+ c_\w |A|L(1+ M^\frac{1}{2}) \sum_{l = 1}^{n-2} M^{-\frac{n-l}{2}} \bigg( \sup_{y \in \X} \frac{\E \big[ \| v(y) - V_l^0(y) \|^2_\infty \big]}{|\w(y)|^2} \bigg)^\frac{1}{2} \nonumber \\
&= |A| M^{- \frac{n}{2}} \big( c_f + c_vLM^\frac{1}{2}\big) + c_\w L(1 + |A|M^{- \frac{1}{2}})F_{n-1} + c_\w |A|L(1+ M^\frac{1}{2}) \sum_{l = 1}^{n-2} M^{-\frac{n-l}{2}} F_l.
\end{align}
This implies for all $n \in \N \cap [2, \infty)$ that
\begin{align} \label{mlp_approx_full_error_gronwall_a_n}
a_n &\leq |A|(c_f + c_vLM^\frac{1}{2}) + c_\w M^\frac{1}{2}L(1 + |A|M^{-\frac{1}{2}}) M^\frac{n-1}{2}F_{n-1} + c_\w |A|L(1+ M^\frac{1}{2})\sum_{l = 1}^{n-2}M^\frac{l}{2}F_l \nonumber \\
&= |A|(c_f + c_vLM^\frac{1}{2}) + c_\w L(M^\frac{1}{2} + |A|)a_{n-1} + c_\w |A|L(1+M^\frac{1}{2})\sum_{l= 1}^{n-2} a_l.
\end{align}
Moreover it holds that 
\begin{align} \label{mlp_approx_full_error_a_0}
a_0 = F_0 = \bigg(\sup_{y \in \X} \frac{\E \big[ \|v(y) \|^2_\infty \big]}{|\w(y)|^2}  \bigg)^\frac{1}{2} = \sup_{(y,a) \in \X \times A} \frac{\left| \dkl{v}{y}{a} \right|}{|\w(y)|}.
\end{align}
The assumption that for all $(x,a) \in \X \times A$ it holds that $ \dkl{v}{x}{a} = \E \big[ f(X^{0,x,a}, v(X^{0,x,a})) \big]$, Jensen's inequality, the triangle inequality, the assumption that for all $x \in \X$, $r,s \in \R^A$ it holds that $|f(x,r) - f(x,s)| \leq L \max_{a \in A} |r(a) - s(a)|$, and the assumption that for all $(x,a) \in \X \times A$ it holds that $\big( \E \big[ |f(X^{0,x,a}, 0|^2 \big] \big)^\frac{1}{2} \leq c_f \w(x)$ and $\big( \E \big[ \| v(X^{0,x,a}) \|_\infty^2 \big)^\frac{1}{2} \leq c_v \w(x)$ prove for all $(x,a) \in \X \times A$ that
\begin{align}
\frac{\big| \dkl{v}{x}{a} \big|}{|\w(x)|} &= \frac{1}{\w(x)} \big( \big| \E \big[ f(X^{0,x,a}, v(X^{0,x,a})) \big] \big|^2 \big)^\frac{1}{2} \leq \frac{1}{\w(x)} \big( \E \big[ \big| f(X^{0,x,a}, v(X^{0,x,a}))\big|^2 \big] \big)^\frac{1}{2} \nonumber \\
&\leq \frac{1}{\w(x)} \Big[ L \big(\E \big[ \sup_{b \in A} |v(X^{0,x,a})(b)|^2 \big] \big)^\frac{1}{2} + \big( \E \big[ |f(X^{0,x,a}, 0)|^2 \big] \big)^\frac{1}{2} \Big] \leq c_vL + c_f
\end{align}
This and (\ref{mlp_approx_full_error_a_0}) ensure that $a_0 \leq c_vL + c_f$. Note that (\ref{mlp_approx_full_error_eq_2}) implies that $a_1 \leq M^\frac{1}{2} \big( c_\w La_0 + |A|M^{-\frac{1}{2}}c_f \big) \leq c_\w L M^\frac{1}{2}(c_vL + c_f) + |A|c_f$. Moreover (\ref{mlp_approx_full_error_gronwall_a_n}) ensures $a_2 \leq |A|(c_f + c_vLM^\frac{1}{2}) + c_\w L(M^\frac{1}{2} + |A|)a_1$. Let $(\xi_n)_{n \in \N_0}$, $(b_n)_{n \in \N_0} \subseteq \R$ satisfy for all $n \in \N_0$ that
\begin{align}\label{mlp_approx_full_error_def_gron}
b_0 &= \xi_0 = \max \Big\{ \frac{c_v}{c_\w}, c_vL + c_f, \frac{|A|c_f}{c_\w |A|L + 1} \Big\}, \quad b_1 = |A|c_f - (c_\w |A|L + 1)b_0,  \nonumber\\
b_2 &= |A|LM^\frac{1}{2}(c_v - c_\w b_0), \quad b_{n+3} = 0,\quad \xi_1 = c_\w LM^\frac{1}{2}b_0 + |A|c_f, \nonumber\\
\xi_{n+2} &= |A|(c_f + c_vLM^\frac{1}{2}) + c_\w L(M^\frac{1}{2} + |A|)\xi_{n+1} + c_\w |A|L(1 + M^\frac{1}{2})\sum_{l = 1}^{n}\xi_l. 
\end{align}
Combining (\ref{mlp_approx_full_error_eq_2}), (\ref{mlp_approx_full_error_gronwall_a_n}), (\ref{mlp_approx_full_error_def_gron}), and induction establishes for all $n \in \N_0$ that $a_n \leq \xi_n$. Observe that (\ref{mlp_approx_full_error_def_gron}) ensures that for all $n \in \N_0$ it holds that
\begin{align}\label{mlp_approx_full_error_two_step}
\xi_0 &= b_0, \quad \xi_1 = b_1 + (c_\w L(M^\frac{1}{2} + |A|) + 1)\xi_0, \nonumber\\
\xi_{n + 2} &= b_{n+2} + (c_\w L(M^\frac{1}{2} + |A|) + 1)\xi_{n+1} + c_\w LM^\frac{1}{2}(|A| - 1)\xi_n.
\end{align} 
Let $x_1, x_2 \in \R$ satisfy 
\begin{align} \label{mlp_approx_full_error_def_x_i} 
x_1 &= \frac{c_\w L (M^\frac{1}{2} + |A|) + 1 - \sqrt{\big( c_\w L (M^\frac{1}{2} + |A|) + 1 \big)^2 + 4c_\w L M^\frac{1}{2} (|A| - 1)  }  }{2}, \nonumber \\
x_2 &= \frac{c_\w L (M^\frac{1}{2} + |A|) + 1 + \sqrt{\big( c_\w L (M^\frac{1}{2} + |A|) + 1 \big)^2 + 4c_\w L M^\frac{1}{2} (|A| - 1)  }  }{2}.
\end{align}
Note that for all $i \in \{ 1,2 \}$ it holds that $x_i^2 = \big(c_\w L (M^\frac{1}{2} + |A|) + 1\big) x_i + c_\w L M^\frac{1}{2}(|A| - 1)$. Moreover note that the assumption that $A$ is nonempty yields that $c_\w L(M^\frac{1}{2} + |A|) + 1 \geq 1$ and $c_\w LM^\frac{1}{2}(|A| - 1) \geq 0$. This ensures that $x_2 > 0 \geq x_1$. This and (\ref{mlp_approx_full_error_def_x_i}) imply that $|x_2| \geq |x_1|$. 
Hence it holds that
\begin{align}
\frac{|x_2|}{|x_2 - x_1|} \leq 1 \quad \text{ and } \quad \frac{|x_1|}{|x_2 - x_1|} \leq \frac{1}{2}.
\end{align}
This, the fact that $|x_1| \leq |x_2|$, and the triangle inequality prove that for all $n \in \N_0$ it holds that
\begin{align} \label{mlp_approx_full_error_x_n_ineq} 
\frac{x_2^{n+1} - x_1^{n+1}}{x_2 - x_1} &\leq \frac{|x_2|^{n + 1}}{x_2 - x_1} + \frac{|x_1|^{n+1}}{x_2- x_1} \leq \frac{x_2}{x_2 - x_1} x_2^n + \frac{|x_1|}{x_2 - x_1} x_2^n \leq \frac{3}{2}x_2^n.
\end{align}
The discrete Gronwall-type two-step recursion in \cite[Lemma 2.1]{HKN2021} (applied with $\kappa \leftarrow c_\w L(M^\frac{1}{2} + |A|) + 1$, $\lambda \leftarrow c_\w LM^\frac{1}{2}(|A| - 1)$, $(a_k)_{k \in \N_0} \leftarrow (\xi_k)_{k \in \N_0}$, $(b_k)_{k \in \N_0} \leftarrow (b_k)_{k \in \N_0}$, $x_{1/2} \leftarrow x_{1/2}$ in the notation of \cite[Lemma 2.1]{HKN2021}) demonstrates that for all $n \in \N_0$ it holds that
\begin{align}
\xi_n = \frac{1}{x_2 - x_1} \Big( b_0 (x_2^{n+1} - x_1^{n+1}) + b_1(x_2^n - x_1^n) +  b_2(x_2^{ \max\{ n-1, 0\} } - x_1^{ \max\{ n-1, 0 \} } ) \Big).
\end{align}
Note that (\ref{mlp_approx_full_error_def_gron}) ensures that $b_1 \leq 0$ and $ b_2 \leq 0$. This, the fact that $x_2 > x_1$, and the fact that for all $n \in \N_0$ it holds that $a_n \leq \xi_n$ proves that for all $n \in \N_0$ it holds that
\begin{align} \label{mlp_approx_full_error_a_n_post_gron}
a_n \leq b_0\frac{x_2^{n+1} - x_1^{n+1}}{x_2 - x_1}.
\end{align}
Combining this and (\ref{mlp_approx_full_error_x_n_ineq}) implies that for all $n \in \N_0$ it holds that
\begin{align} \label{mlp_approx_full_error_f_n_post_gron}
F_n \leq \frac{3}{2}b_0 \big( M^{-\frac{1}{2}} x_2 \big)^n.
\end{align}
Furthermore, observe that
\begin{align}
M^{-\frac{1}{2}}x_2 &= M^{-\frac{1}{2}} \Bigg( \frac{ c_\w L(M^\frac{1}{2} + |A|) + 1 + \sqrt{\big( c_\w L(M^\frac{1}{2} + |A|) + 1 \big)^2 + 4c_\w LM^\frac{1}{2}(|A| - 1) }   }{2} \Bigg) \nonumber \\
&= \frac{ c_\w L(1 + |A|M^{-\frac{1}{2}}) + M^{-\frac{1}{2}} + \sqrt{\big( c_\w L(1 + |A|M^{-\frac{1}{2}}) + M^{-\frac{1}{2}} \big)^2 + 4c_\w LM^{-\frac{1}{2}}(|A| - 1)}   }{2}
\end{align}
This, (\ref{mlp_approx_full_error_def_fn_an}), and (\ref{mlp_approx_full_error_f_n_post_gron}) prove that for all $x \in \X$, $n \in \N$ it holds that
\begin{align}
&\bigg( \frac{\E \big[ \| v(x) - V_n^\theta(x) \|_\infty^2 \big]}{|\w(x)|^2} \bigg)^\frac{1}{2} \\
&\leq c \Bigg( \frac{c_\w L(1 + |A|M^{-\frac{1}{2}}) + M^{-\frac{1}{2}} + \sqrt{\big( c_\w L(1 + |A|M^{-\frac{1}{2}}) + M^{-\frac{1}{2}}\big)^2 + 4c_\w LM^{-\frac{1}{2}}(|A| - 1)}  }{2}  \Bigg)^n. \nonumber
\end{align}
The proof of \cref{mlp_approx_full_error} is thus completed.
\end{mproof}

\begin{corollary} \label{mlp_approx_convergence} %(Convergence)(c.f. cor 3.9 in stopping)\\
Assume \cref{mlp_setting}, assume $A$ is finite, let $c_f, c_v, c_\w, L \in [0, \infty)$ with $c_\w L < 1$, let $M > (1 + c_\w L(2|A| - 1))^2(1-c_\w L)^{-2}$, let $\w \colon \X \rightarrow (0, \infty)$ be $\mathcal{X}/\mathcal{B}((0, \infty))$-measurable, let $v \colon \X \rightarrow \R^A$ be $\mathcal{X} / \A$-measurable, assume for all $(x,a) \in \X \times A$, $r,s \in \R^A$ that $|f(x,r) - f(x,s)| \leq L \max_{b \in A} |r(b) - s(b)|$, $\big( \E \big[ |f(X^{0,x,a}, 0)|^2 \big] \big)^\frac{1}{2} \leq c_f \w(x)$, $\big( \E \big[ \|v(X^{0,x,a}) \|_\infty^2 \big] \big)^\frac{1}{2} \leq c_v \w(x)$, $\big( \E \big[ \big| \w(X^{0,x,a}) \big|^2 \big] \big)^\frac{1}{2} \leq c_\w \w(x)$, and $\dkl{v}{x}{a} = \E\big[ f\big( X^{0,x,a}, v(X^{0,x,a}) \big) \big]$. Then it holds that 
\begin{align}
\limsup_{n \rightarrow \infty}  \bigg[ \sup_{x \in \X} \bigg( \frac{ \E \big[ \| v(x) - V_n^0(x) \|^2_\infty \big] }{|\w(x)|^2} \bigg)^\frac{1}{2} \bigg] = 0.
\end{align}
\end{corollary}

\begin{mproof}{\cref{mlp_approx_convergence}}
Let $c \in [0, \infty)$ satisfy $c = \frac{3}{2} \max \big\{ \frac{c_v}{c_\w}, c_vL + c_f, \frac{|A|c_f}{c_\w |A|L + 1} \big\} $. \cref{mlp_approx_full_error} establishes that for all $n \in \N$, $x \in \X$ it holds that
\begin{align} \label{mlp_approx_convergence_base_ineq} 
&\bigg( \frac{\E \big[ \| v(x) - V_n^0(x)\|^2_\infty \big]}{|\w(x)|^2} \bigg)^\frac{1}{2} \\
&\leq c \bigg( \frac{c_\w L(1 + |A|M^{-\frac{1}{2}}) + M^{-\frac{1}{2}} + \sqrt{\big(c_\w L(1 + |A|M^{-\frac{1}{2}}) + M^{-\frac{1}{2}}\big)^2 + 4c_\w LM^{-\frac{1}{2}}(|A| - 1)}  }{2}  \Bigg)^n.  \nonumber
\end{align}   
Observe that the assumption that $M > (1 + c_\w L(2|A| - 1))^2(1- c_\w L)^{-2}$ implies $M^{-\frac{1}{2}} < (1- c_\w L)(1 + c_\w L(2|A| - 1))^{-1}$. This proves that
\begin{align} \label{mlp_approx_convergence_aux_eq_1}
c_\w L(|A| - 1)M^{-\frac{1}{2}} &< \frac{c_\w L(|A| - 1)(1- c_\w L)}{1 + c_\w L(2|A| - 1)} \nonumber\\
&= 1 - \Big( c_\w L + (1 - c_\w L) - \frac{c_\w L (|A| - 1) (1- c_\w L)}{1 + c_\w L(2 |A| - 1)} \Big) \nonumber \\
&= 1 - \Big( c_\w L + \big( 1 + c_\w L (2 |A| - 1) - c_\w L (|A| - 1) \big) \frac{1 - c_\w L}{1 + c_\w L (2|A| - 1)}  \Big)  \nonumber \\
&= 1 - \Big( c_\w L + \big( 1 + c_\w |A|L  \big) \frac{1 - c_\w L}{1 + c_\w L (2|A| - 1)}  \Big)  \nonumber \\
&< 1 - \big( c_\w L + (1 + c_\w |A|L)M^{-\frac{1}{2}} \big) = 1 - \big(c_\w L(1 + |A|M^{-\frac{1}{2}}) + M^{-\frac{1}{2}} \big).
\end{align}
The assumption that $A$ is nonempty ensures that $|A| \leq 2|A| - 1$. This and the assumption that $M > (1 + c_\w L(2|A| - 1))^2(1-c_\w L)^{-2}$ yield $M > (1 + c_\w L|A|)^2(1-c_\w L)^{-2}$. Therefore it holds that $c_\w L(1 + |A|M^{-\frac{1}{2}}) + M^{-\frac{1}{2}} = c_\w L + (1 + c_\w |A|L)M^{-\frac{1}{2}} < 1 < 2$. This and (\ref{mlp_approx_convergence_aux_eq_1}) establish that
\begin{align} \label{mlp_approx_convergence_aux_eq_2}
&\hspace{-1.5cm}\tfrac{c_\w L(1 + |A|M^{-\frac{1}{2}}) + M^{-\frac{1}{2}} + \sqrt{\big( c_\w L(1 + |A|M^{-\frac{1}{2}}) + M^{-\frac{1}{2}}\big)^2 + 4c_\w LM^{-\frac{1}{2}}(|A| - 1)}  }{2} \nonumber\\
&< \tfrac{c_\w L(1 + |A|M^{-\frac{1}{2}}) + M^{-\frac{1}{2}} + \sqrt{\big(c_\w L(1 + |A|M^{-\frac{1}{2}}) + M^{-\frac{1}{2}}\big)^2 + 4 - 4\big(c_\w L(1 + |A|M^{-\frac{1}{2}}) + M^{-\frac{1}{2}} \big) } }{2} \nonumber\\
&= \tfrac{c_\w L(1 + |A|M^{-\frac{1}{2}}) + M^{-\frac{1}{2}} + \sqrt{\big(2 - c_\w L(1 + |A|M^{-\frac{1}{2}}) + M^{-\frac{1}{2}}\big)^2 } }{2} = 1.
\end{align}
Combining (\ref{mlp_approx_convergence_base_ineq}) and (\ref{mlp_approx_convergence_aux_eq_2}) implies that 
\begin{align}
\limsup_{n \rightarrow \infty} \bigg[ \sup_{x \in \X} \bigg( \frac{ \E \big[ \| v(x) - V_n^0(x) \|^2_\infty \big] }{|\w(x)|^2} \bigg)^\frac{1}{2} \bigg] = 0.
\end{align}
The proof of \cref{mlp_approx_convergence} is thus completed.
\end{mproof}

\section{Computational complexity analysis for MLFP approximations}\label{sec:comp_comp}

\subsection{MLFP approximations for functional fixed-point equations}

\begin{theo} \label{mlp_approx_main_theo} %(cf. Thm 3.10 in stopping )\\
Let $M \in \N$, let $\Theta = \bigcup_{n \in \N} \Z^n$, let $A$ be a finite nonempty set, let $\kappa \in [0,\infty)$, let $\Dm$ be a nonempty set, let $\lambda_\dm, L_\dm, \mathfrak{R}_\dm \in [0,\infty)$, $\dm \in \Dm$, let $(\Omega, \mathcal{F}, \mathbb{P})$ be a probability space, let $(\X_\dm, \mathcal{X}_\dm)$, $\dm \in \Dm$, be nonempty measurable spaces, for every $\dm \in \Dm$ let $\w_\dm \colon \X_\dm \rightarrow (0,\infty)$ be $\mathcal{X}_\dm/\mathcal{B}((0,\infty))$-measurable, for every $\dm \in \Dm$ let $f_\dm \colon \X_\dm \times \R^A \rightarrow \R$, be $(\mathcal{X}_\dm \otimes (\bigotimes_{a \in A}\mathcal{B}(\R)) )/\mathcal{B}(\R)$-measurable, for every $\dm \in \Dm$ let $(\mathcal{F}_\dm^\theta)_{\theta \in \Theta}$ be independent sub-$\sigma$-algebras of $\mathcal{F}$, for every $\dm \in \Dm$ let $X_\dm^\theta \colon \X_\dm \times \Omega \rightarrow \X_\dm^A$, $\theta \in \Theta$, be i.i.d.\ random fields which satisfy for all $\dm \in \Dm$, $\theta \in \Theta$ that $X_\dm^\theta$ is $(\mathcal{X}_\dm \otimes \mathcal{F}_\dm^\theta)/(\bigotimes_{a \in A} \mathcal{X}_\dm)$-measurable, assume for all $\dm \in \Dm$, $(x,a) \in \X_\dm \times A$, $r,s \in \R^A$ that $|f_\dm(x,r) - f_\dm(x,s)| \leq L_\dm \max_{b \in A} |r(b) - s(b)|$, $ \left(\E\big[ |f_\dm(X_\dm^{0,x,a}, 0)|^2 \big]\right)^\frac{1}{2} \leq \kappa \w_\dm(x)$, $\left(\E \big[ | \w_\dm(X_\dm^{0,x,a})|^2 \big]\right)^\frac{1}{2} \leq \lambda_\dm \w_\dm(x)$, $\sup_{u \in \Dm} \lambda_u L_u < 1$, 
assume $M > \frac{( 1 + (\sup_{\dm \in \Dm}  \lambda_\dm L_\dm)(2|A| - 1) )^2}{ (1 - (\sup_{\dm \in \Dm} \lambda_\dm L_\dm))^2 }$, for every $\dm \in \Dm$ let $V_{n,\dm}^\theta \colon \X_\dm \times \Omega \rightarrow \R^A$, $n \in \N_0$, $\theta \in \Theta$, satisfy for all $n \in \N_0$, $(x,a) \in \X_\dm \times A$, $\theta \in \Theta$ that
\begin{align} \label{mlp_approx_main_theo_def_v_n}
\dkl{V_{n, \dm}^\theta}{x}{a} = \sum_{l = 0}^{n-1} \frac{1}{M^{n-l}} \sum_{i = 1}^{M^{n-l}} &f_\dm\big( X_\dm^{(\theta, l, i),x ,a}, V_{l, \dm}^{(\theta, l, i)}(X_\dm^{(\theta, l, i), x, a})  \big) \nonumber \\
&-\mathbbm{1}_{\N}(l) f_\dm \big( X_\dm^{(\theta, l, i), x, a}, V_{\max\{ l-1, 0\}, \dm}^{(\theta, -l, i)}(X_\dm^{(\theta, l, i), x, a}) \big),
\end{align}
and let $\mathfrak{C}_{n,\dm} \in [0, \infty)$, $n \in \N_0$, $\dm \in \Dm$, satisfy for all $n \in \N_0$, $\dm \in \Dm$ that
\begin{align} \label{mlp_approx_main_theo_comp_effort}
\mathfrak{C}_{n, \dm} \leq \sum_{l = 0}^{n-1} M^{n-l}\big( \mathfrak{R}_\dm + \mathfrak{C}_{l, \dm} + \mathbbm{1}_{\N}(l) \mathfrak{C}_{\max \{l-1, 0\}, \dm } \big).
\end{align}
Then the following holds: 
\begin{itemize}
\item[$(i)$] For every $\dm \in \Dm$ there exists a unique function $v_\dm \colon \X_\dm \rightarrow \R^A$ which is $\mathcal{X}_\dm / (\bigotimes_{a \in A} \mathcal{B}(\R))$-measurable \hspace{0.05cm} and \hspace{0.05cm} satisfies \hspace{0.05cm} for \hspace{0.05cm} all \hspace{0.2cm} $(x,a) \in \X_\dm \times A$ \hspace{0.2cm} that \hspace{0.28cm} $\sup_{y \in \X_\dm} \frac{\| v_\dm(y) \|_\infty }{\w_\dm (y)} < \infty$, $\E \big[ |f_\dm(X_\dm^{0,x,a}, v_\dm(X_\dm^{0,x,a}))| \big] < \infty$, and 
\begin{align}\label{eq:gen_func_equation}
\dkl{v_\dm}{x}{a} = \E \big[ f_\dm(X_\dm^{0,x,a}, v_\dm(X_\dm^{0,x,a})) \big].
\end{align}
\item[$(ii)$] There exist $N \colon (0, 1] \rightarrow \N$ and $c \in \R$ such that for all $\dm \in \Dm$, $\varepsilon \in (0,1]$ it holds that $\mathfrak{C}_{N_\varepsilon , \dm} \leq c \mathfrak{R}_\dm \varepsilon^{-c}$ and
\begin{align} \label{mlp_approx_main_theo_error} 
\sup_{x \in \X_\dm} \bigg( \frac{\E \big[ \| v_\dm(x) - V_{N_\varepsilon, \dm}^0(x) \|^2_\infty \big]}{|\w_\dm(x)|^2} \bigg)^\frac{1}{2} \leq \varepsilon.
\end{align}
\end{itemize}
\end{theo}

\begin{mproof}{\cref{mlp_approx_main_theo}}
Let $\alpha, \beta, \gamma \in \R \cup \{\infty \}$ satisfy 
\begin{align} \label{mlp_approx_main_theo_constants} 
\alpha &= \frac{1}{2} \sup_{\dm \in \Dm} \Big[ \lambda_\dm L_\dm(1 + |A|M^{-\frac{1}{2}}) + M^{-\frac{1}{2}} + \sqrt{\big( \lambda_\dm L_\dm(1 + |A|M^{-\frac{1}{2}}) + M^{-\frac{1}{2}} \big)^2 + 4M^{-\frac{1}{2}} \lambda_\dm L_\dm(|A| - 1) } \Big], \nonumber \\
\beta &= \frac{\ln(3M)}{\ln(\alpha^{-1})}, \quad \text{and } \quad \gamma = \frac{3}{2} \sup_{\dm \in \Dm} \left( \max \Big\{ \frac{\kappa}{1 - \lambda_\dm  L_\dm}, \frac{\kappa \lambda_\dm L_\dm}{1 - \lambda_\dm L_\dm} + \kappa, \frac{|A|\kappa}{\lambda_\dm |A|L_\dm + 1}  \Big\}\right).
\end{align}
Let $N \colon (0,1] \rightarrow \N \cup \{ \infty \}$ satisfy for all $\varepsilon \in (0,1]$ that $N_\varepsilon = \min (\{ n \in \N : \gamma \alpha^n \leq \varepsilon  \} \cup \{ \infty \})$. The assumption that $\sup_{\dm \in \Dm} \lambda_\dm L_\dm < 1$ ensures that $\gamma \in [0, \infty)$. The triangle inequality yields
\begin{align} \label{mlp_approx_main_theo_a_tri_ineq} 
\alpha &= \frac{1}{2} \sup_{\dm \in \Dm} \Big[ \lambda_\dm L_\dm(1 + |A|M^{-\frac{1}{2}}) + M^{-\frac{1}{2}} + \sqrt{\big( \lambda_\dm L_\dm(1 + |A|M^{-\frac{1}{2}}) + M^{-\frac{1}{2}} \big)^2 + 4M^{-\frac{1}{2}} \lambda_\dm L_\dm(|A| - 1) } \Big], \nonumber \\
&\leq \frac{1}{2} \bigg[ (\sup_{\dm \in \Dm} \lambda_\dm L_\dm )(1 + |A|M^{-\frac{1}{2}}) + M^{-\frac{1}{2}} \nonumber\\
&\hspace{1cm}+ \sqrt{\big( ( \sup_{\dm \in \Dm} \lambda_\dm L_\dm )(1 + |A|M^{-\frac{1}{2}}) + M^{-\frac{1}{2}} \big)^2 + 4M^{- \frac{1}{2}} (\sup_{\dm \in \Dm}\lambda_\dm L_\dm)(|A| - 1) } \bigg]. 
\end{align}
The assumption that $M^{-\frac{1}{2}} < \frac{1 - (\sup_{\dm \in \Dm} \lambda_\dm L_\dm)}{1 + (\sup_{\dm \in \Dm} \lambda_\dm L_\dm)(2|A| - 1)}$ ensures that
\begin{align} \label{mlp_approx_main_theo_help_ineq}
1 - \big( (\sup_{\dm \in \Dm} \lambda_\dm L_\dm) (1 + |A|M^{-\frac{1}{2}}) + M^{-\frac{1}{2}} \big) &= 1 - (\sup_{\dm \in \Dm} \lambda_\dm L_\dm) - (1 + |A| (\sup_{\dm \in \Dm} \lambda_\dm L_\dm )) M^{-\frac{1}{2}}  \nonumber\\
&\hspace{-5.5cm}> (1 -  (\sup_{\dm \in \Dm} \lambda_\dm L_\dm)) - (1 + |A| (\sup_{\dm \in \Dm} \lambda_\dm L_\dm )) \frac{1 - (\sup_{\dm \in \Dm} \lambda_\dm L_\dm)}{1 + (\sup_{\dm \in \Dm} \lambda_\dm L_\dm)(2|A| - 1)}  \nonumber\\
&\hspace{-5.5cm}= \Big( 1 + (\sup_{\dm \in \Dm} \lambda_\dm L_\dm)(2 |A| - 1) - \big( |A|(\sup_{\dm \in \Dm} \lambda_\dm L_\dm) + 1 \big) \Big) \frac{1 - (\sup_{\dm \in \Dm} \lambda_\dm L_\dm)}{1 + (\sup_{\dm \in \Dm} \lambda_\dm L_\dm)(2|A| - 1)} \nonumber \\
&\hspace{-5.5cm}= \Big( (\sup_{\dm \in \Dm} \lambda_\dm L_\dm) (|A| - 1)\Big) \frac{1 - (\sup_{\dm \in \Dm} \lambda_\dm L_\dm)}{1 + (\sup_{\dm \in \Dm} \lambda_\dm L_\dm)(2|A| - 1)} \nonumber\\
&\hspace{-5.5cm}> (\sup_{\dm \in \Dm} \lambda_\dm L_\dm) (|A| - 1) M^{-\frac{1}{2}}.
\end{align}
The assumption that $A$ is nonempty implies that $|A| \leq 2|A| - 1$. Hence in holds that $M > \frac{(1 + (\sup_{\dm \in \Dm} \lambda_\dm L_\dm)(2|A| - 1) )^2}{(1 - (\sup_{\dm \in \Dm}\lambda_\dm L_\dm))^2 } \geq \frac{(1 + |A|(\sup_{t \in T}\lambda_\dm L_\dm)^2}{(1 - (\sup_{\dm \in \Dm} \lambda_\dm L_\dm))^2}$. This establishes that
\begin{align}
(\sup_{\dm \in \Dm} \lambda_\dm L_\dm) + (|A|(\sup_{\dm \in \Dm} \lambda_\dm L_\dm) + 1) M^{-\frac{1}{2}} < 1.
\end{align} 
Combining this, (\ref{mlp_approx_main_theo_a_tri_ineq}), and (\ref{mlp_approx_main_theo_help_ineq}) demonstrates that
\begin{align}
\alpha &\leq \frac{1}{2} \bigg[ (\sup_{\dm \in \Dm} \lambda_\dm L_\dm) (1 + |A|M^{-\frac{1}{2}}) + M^{-\frac{1}{2}} \nonumber\\
&\hspace{0.8cm}+ \sqrt{\big( (\sup_{\dm \in \Dm} \lambda_\dm L_\dm) (1 + |A|M^{-\frac{1}{2}}) + M^{-\frac{1}{2}} \big)^2 + 4M^{- \frac{1}{2}} (\sup_{\dm \in \Dm}\lambda_\dm L_\dm)( |A| - 1) } \bigg] \nonumber\\
&< \frac{1}{2} \bigg[ (\sup_{\dm \in \Dm} \lambda_\dm L_\dm) (1 + |A|M^{-\frac{1}{2}}) + M^{-\frac{1}{2}} \nonumber\\
&\hspace{0.8cm}+ \sqrt{\big( (\sup_{\dm \in \Dm} \lambda_\dm L_\dm) (1 + |A|M^{-\frac{1}{2}}) + M^{-\frac{1}{2}} \big)^2 + 4\big( 1 - (\sup_{\dm \in \Dm} \big(\lambda_\dm L_\dm)(1 + |A|M^{-\frac{1}{2}}) + M^{-\frac{1}{2}} \big)  \big) } \bigg] \nonumber\\
&= \frac{1}{2} \bigg[ (\sup_{\dm \in \Dm} \lambda_\dm L_\dm) (1 + |A|M^{-\frac{1}{2}}) + M^{-\frac{1}{2}} + \sqrt{\big(2 - \big( (\sup_{\dm \in \Dm} \lambda_\dm L_\dm) (1 + |A|M^{-\frac{1}{2}}) + M^{-\frac{1}{2}} \big) \big)^2 } \bigg] \nonumber\\
&= 1.
\end{align}
This and (\ref{mlp_approx_main_theo_constants}) ensure that $\alpha \in [M^{-\frac{1}{2}}, 1)$. The assumption that $M > \frac{(1 + (\sup_{\dm \in \Dm} \lambda_\dm L_\dm) (2 |A| - 1) )^2}{(1- (\sup_{\dm \in \Dm}\lambda_\dm L_\dm))^2}$ implies that $M > 1$. Hence it holds that $[M^{- \frac{1}{2}}, 1) \neq \emptyset$. The fact that $\alpha \in (0,1)$ demonstrates that for all $\varepsilon \in (0,1]$ it holds that $N_\varepsilon \in \N$. The fact that $\alpha \in [M^{-\frac{1}{2}}, 1)$ and the fact that $M > 1$ ensure that $\beta \in (2, \infty)$. \cref{ex_unique_sol_sto_fp_eq_ew} yields that for all $\dm \in \Dm$ there exists a unique function $v_\dm \colon \X_\dm \rightarrow \R^A$ which is $\mathcal{X}_\dm/ (\bigotimes_{a \in A} \mathcal{B}(\R))$-measurable and satisfies for all $\dm \in \Dm$, $(x,a) \in \X_\dm \times A$ that $\sup_{y \in \X_\dm} \frac{\| v_\dm(y)\|_\infty}{\w_\dm(y)} < \infty$, $\E \big[ |f_\dm(X_\dm^{0,x,a} v_\dm(X_\dm^{0,x,a}))| \big] < \infty$, and $v_\dm(x,a) = \E \big[ f_\dm(X_\dm^{0,x,a}, v_\dm(X_\dm^{0,x,a})) \big]$. This proves item $(i)$. \cref{sol_sto_fp_eq_estimates} implies for all $\dm \in \Dm$, $(x,a) \in \X_\dm \times A$ that $ \big( \E \big[ \| v_\dm(X_\dm^{0,x,a}) \|_\infty^2\big]\big)^\frac{1}{2} \leq \frac{\kappa \lambda_\dm}{1 - \lambda_\dm L_\dm} \w_\dm(x)$.
\cref{mlp_approx_full_error} establishes for all $n \in \N$, $\dm\in \Dm$, $x \in \X_\dm$ that
\begin{align}
&\bigg(\frac{\E \big[ \| v_\dm(X_\dm^{0,x,a}) - V_{n,\dm}^0(X_\dm^{0,x,a}) \|_\infty^2 \big]}{|\w_\dm(x)|^2} \bigg)^\frac{1}{2} \\
&\hspace{1cm}\leq \gamma \left( \tfrac{ \lambda_\dm L_\dm (1 + |A|M^{-\frac{1}{2}}) + M^{-\frac{1}{2}} + \sqrt{(\lambda_\dm L_\dm(1 + |A|M^{-\frac{1}{2}}) + M^{-\frac{1}{2}})^2 + 4\lambda_\dm L_\dm M^{-\frac{1}{2}} (|A| - 1)  }}{2} \right)^n \hspace{-0.3cm}.  \nonumber
\end{align}
This demonstrates for all $n \in \N$, $\dm \in \Dm$ that
\begin{align} \label{mlp_approx_main_theo_error_w_a} 
\sup_{x \in \X_\dm} \bigg( \frac{\E \big[ \| v_\dm(x) - V_{n,\dm}^0(x) \|_\infty^2 \big]}{|\w_\dm(x)|^2} \bigg)^\frac{1}{2} \leq \gamma \alpha^n.
\end{align}
Note that \cite[Lemma 3.14]{beck2020overcoming} (applied with $M \leftarrow M$, $\alpha \leftarrow \mathfrak{R}_\dm + \mathfrak{C}_{0,\dm}$, $\beta \leftarrow \mathfrak{R}_\dm$, $(C_n) \leftarrow (\mathfrak{C}_{n,\dm})$ for $\dm \in \Dm$ in the notation of \cite[Lemma 3.14]{beck2020overcoming}) proves that for all $n \in \N$, $\dm \in \Dm$ it holds that
\begin{align}
\mathfrak{C}_{n, \dm} \leq \bigg( \frac{\mathfrak{R}_\dm + \mathfrak{C}_{0,\dm} + \mathfrak{R}_\dm + \mathfrak{C}_{0,\dm} }{2} \bigg) (3M)^n = \mathfrak{R}_\dm (3M)^n.
\end{align} 
This, (\ref{mlp_approx_main_theo_error_w_a}), and \cite[Lemma 3.15]{beck2020overcoming} (applied with $m \leftarrow 1$, $\alpha \leftarrow \alpha$, $\beta \leftarrow 3M$, $\kappa_1 \leftarrow \gamma$, $\kappa_2 \leftarrow \mathfrak{R}_\dm$, $N \leftarrow N$ for $\dm \in \Dm$ in the notation of \cite[Lemma 3.15]{beck2020overcoming}) implies for all $\dm \in \Dm$, $\varepsilon \in (0,1]$ that 
\begin{align}
\sup_{x \in \X_\dm} \bigg( \frac{\E \big[ \| v_\dm(x) - V_{n, \dm}^0(x) \|^2_\infty \big]}{|\w_\dm(x)|^2} \bigg)^\frac{1}{2} \leq \varepsilon \quad \text{and} \quad \mathfrak{C}_{N_\varepsilon, \dm} \leq 3M\mathfrak{R}_\dm \max \{ 1, \gamma \}^\beta \frac{1}{\varepsilon^\beta}.
\end{align}
Let $c = \max\{ \beta, 3M \max\{1,\gamma \}^\beta \}$. Hence it holds for all $\dm \in \Dm$, $\varepsilon \in (0,1]$ that $\mathfrak{C}_{N_\varepsilon, \dm} \leq c \mathfrak{R}_\dm \varepsilon^{-c}$. This completes the proof of \cref{mlp_approx_main_theo}.
\end{mproof}

\subsection{MLFP approximations for Bellman equations of optimal control problems}

\begin{corollary} \label{mlp_approx_cor_ctrl}
    Let $M \in \N$, $\kappa \in [0,\infty)$, $\Theta = \bigcup_{n \in \N} \Z^n$, let $A$ be a finite nonempty set, let $\Dm$ be a nonempty set, let $(\Omega, \mathcal{F}, \mathbb{P})$ be a probability space, 
    let $\lambda_{\dm}, \delta_\dm, \mathfrak{R}_\dm \in [0,\infty)$, $\dm \in \Dm$,
    let $(\X_\dm, \mathcal{X}_\dm)$, $\dm \in \Dm$, be nonempty measurable spaces,
    for every $\dm \in \Dm$ let $\w_\dm \colon \X_\dm \rightarrow (0, \infty)$ be $\mathcal{X}_\dm/\mathcal{B}((0,\infty))$-measurable, 
    for every $\dm \in \Dm$ let $g_\dm \colon \X_\dm \times A \rightarrow \R$ be $(\mathcal{X}_\dm \otimes 2^A )/\mathcal{B}(\R)$-measurable, 
    for every $\dm \in \Dm$ let $(\mathcal{F}^\theta_\dm)_{\theta \in \Theta}$ be independent sub-$\sigma$-algebras of $\mathcal{F}$, for every $\dm \in \Dm$ let $X^\theta_\dm = \big( X^{\theta, x, a}_\dm(\omega) \big)_{x \in \X_\dm, \;a \in A, \; \omega \in \Omega} \colon \X_\dm \times A \times \Omega\rightarrow \X_\dm$, $\theta \in \Theta$, be i.i.d.\ random fields which satisfy for all $\dm \in \Dm$, $\theta \in \Theta$ that $X^\theta_\dm$ is $(\mathcal{X}_\dm \otimes 2^A \otimes \mathcal{F}^\theta_\dm)/\mathcal{X}_\dm$-measurable, assume $M > \frac{( 1 + (\sup_{\dm \in \Dm} \lambda_\dm \delta_\dm)(2|A| - 1) )^2}{ (1 - (\sup_{\dm \in \Dm}\lambda_\dm \delta_\dm))^2 }$, assume for all $\dm \in \Dm$, $(x,a) \in \X_\dm \times A$ that $\max_{b \in A} | g_\dm(x,b) |\leq \kappa \w_\dm(x)$,  $\big( \E \big[ \left|\w_\dm(X_\dm^{0,x,a})\right|^2 \big] \big)^\frac{1}{2} \leq \lambda_\dm \w_\dm(x)$, $\sup_{u \in \Dm} \lambda_u \delta_u < 1$, 
    for every $\dm \in \Dm$ let $Q^\theta_{n, \dm} \colon \X_\dm \times A \times \Omega \rightarrow \R$, $\theta \in \Theta$, $n \in \N_0$, satisfy for all $n \in \N_0$, $(x,a) \in \X_\dm \times A$, $\theta \in \Theta$ that
    \begin{align}
        Q_{n, \dm}^{\theta}(x,a) = g_\dm(x,a) + \sum_{l = 0}^{n-1} \frac{\delta_\dm}{M^{n-l}} \sum_{i = 1}^{M^{n-l}} &\max_{b \in A}\big\{ Q_{l, \dm}^{(\theta, l, i)}(X_\dm^{(\theta, l, i), x,a},b) \big\} \\
        &\hspace{-0cm}-\mathbbm{1}_{\N}(l) \max_{b \in A} \big\{ Q_{\max\{ l-1, 0 \}, \dm}^{(\theta, -l, i)}(X_\dm^{(\theta, l, i), x,a},b) \big\}, \nonumber 
    \end{align}
    and let $\mathfrak{C}_{n,\dm} \in [0, \infty)$, $n \in \N_0$, $\dm \in \Dm$, satisfy for all $n \in \N_0$, $\dm \in \Dm$ that
    \begin{align}\label{eq:comp_effort_cor_bellm}
        \mathfrak{C}_{n, \dm} \leq \sum_{l = 0}^{n-1} M^{n-l}\big( \mathfrak{R}_\dm + \mathfrak{C}_{l, \dm} + \mathbbm{1}_{\N}(l) \mathfrak{C}_{\max \{l-1, 0\}, \dm } \big).
    \end{align}
    Then the following holds:
    \begin{itemize}
        \item[$(i)$] For every $\dm \in \Dm$ there exists a unique function $Q_\dm \colon \X_\dm \times A \rightarrow \R$ which is $(\mathcal{X}_\dm \otimes 2^A )/ \mathcal{B}(\R)$-measurable \hspace{0.1cm} and \hspace{0.1cm} satisfies \hspace{0.1cm} for \hspace{0.1cm} all \hspace{0.2cm} $(x,a) \in \X_\dm \times A$ \hspace{0.1cm} that \hspace{0.1cm} $\sup_{y \in \X_\dm} \frac{\max\limits_{b \in A}\left| Q_\dm(y, b) \right|}{\w_\dm (y)} < \infty$, $\E \big[ \big| \max_{b \in A} Q_\dm(X_\dm^{0,x,a},b)  \big| \big] < \infty$, and 
        \begin{align}
            Q_\dm(x,a) = g_\dm (x,a) + \delta_\dm \E \big[ \max_{b \in A} Q_d(X_\dm^{0,x,a}, b) \big].
        \end{align}
        \item[$(ii)$] There exist $N \colon (0, 1] \rightarrow \N$ and $c \in \R$ such that for all $\dm \in \Dm$, $\varepsilon \in (0,1]$ it holds that $\mathfrak{C}_{N_\varepsilon , \dm} \leq c \mathfrak{R}_\dm \varepsilon^{-c}$ and
        \begin{align}
            \sup_{x \in \X_\dm} \bigg( \frac{\E \big[ \max_{a \in A}| Q_\dm(x,a) - Q_{N_\varepsilon, \dm}^0(x, a) |^2 \big]}{|\w_\dm(x)|^2} \bigg)^\frac{1}{2} \leq \varepsilon.
        \end{align}
    \end{itemize}
\end{corollary}

\begin{mproof}{\cref{mlp_approx_cor_ctrl}}
Note that for every $\dm \in \Dm$ it holds that the function $\X_\dm \times \R^A \ni (x,r) \mapsto \delta_\dm \max_{a \in A} \left\{ g_\dm(x,a) + r(a) \right\} \in \R$ is $(\mathcal{X}_\dm \otimes ( \bigotimes_{a \in A} \mathcal{B}(\R)))/\mathcal{B}(\R)$-measurable. Moreover, for all $\dm \in \Dm$, $x \in \X_\dm$, $r, s \in \R^A$ it holds that 
\begin{align}
\big| \delta_\dm \max_{a \in A} \left\{ g_\dm(x,a) + r(a)  \right\} - \delta_\dm \max_{a \in A} \left\{ g_\dm(x,a) + s(a) \right\}\big| &\leq \delta_\dm \max_{a \in A} \left| r(a) - s(a) \right|.
\end{align}
The assumption that for all $\dm \in \Dm$, $(x,a) \in \X_\dm \times A$ it holds that $\max_{b \in A} | g_\dm(x,b) |\leq \kappa \w_\dm(x)$,  $\big( \E \big[ \left|\w_\dm(X_\dm^{0,x,a})\right|^2 \big] \big)^\frac{1}{2} \leq \lambda_\dm \w_\dm(x)$, and $\sup_{u \in \Dm} \lambda_u \delta_u < 1$ yields for all $\dm \in \Dm$, $(x,a) \in \X_\dm \times A$ that
\begin{align}
\left( \E \left[ \left| \delta_\dm \max_{b \in A} \left\{ g_\dm(X_\dm^{0,x,a}, b) \right\} \right|^2 \right] \right)^\frac{1}{2} \leq \delta_\dm \kappa \left( \E \left[ \left| \w_\dm(X_\dm^{0,x,a}) \right|^2 \right] \right)^\frac{1}{2} \leq \delta_\dm \kappa \lambda_\dm \w_\dm(x) \leq \kappa \w_\dm(x).
\end{align}
For every $\dm \in \Dm$ let $R_{n, \dm}^\theta \colon \X_\dm \times A \rightarrow \R$, $\theta \in \Theta$, $n \in \N_0$, satisfy for all $n \in \N_0$, $(x,a) \in \X_\dm \times A$, $\theta \in \Theta$ that
    \begin{align}\label{eq:mlp_approx_cor_ctrl_R_scheme} 
        R_{n, \dm}^\theta(x,a) &= \sum_{l = 0}^{n-1} \frac{\delta_\dm}{M^{n-l}} \sum_{i = 1}^{M^{n - l}} \max_{b \in A} \left\{ g_\dm (X_\dm^{(\theta, l, i), x, a}, b) + R_{l, \dm}^{(\theta, l, i)}\big( X_\dm^{(\theta, l, i), x, a}, b \big) \right\} \nonumber \\
        &\hspace{3cm}- \mathbbm{1}_{\N}(l) \max_{b \in A} \left\{ g_\dm (X_\dm^{(\theta, l, i), x, a}, b) + R_{\max \{ l-1, 0 \}, \dm}^{(\theta, -l, i)}\big( X_\dm^{(\theta, l, i), x, a}, b \big) \right\}.
    \end{align}
    This and \cref{mlp_approx_main_theo} (applied with $M \leftarrow M$, $\Theta \leftarrow \Theta$, $A \leftarrow A$, $\kappa \leftarrow \kappa$, $(\Omega, \mathcal{F}, \mathbb{P}) \leftarrow (\Omega, \mathcal{F}, \mathbb{P})$, $\Dm \leftarrow \Dm$, $(\lambda_\dm)_{\dm \in \Dm} \leftarrow (\lambda_\dm)_{\dm \in \Dm}$, $(L_\dm)_{\dm \in \Dm} \leftarrow (\delta_\dm)_{\dm \in \Dm}$, $(\mathfrak{R}_\dm)_{\dm \in \Dm} \leftarrow (\mathfrak{R}_\dm)_{\dm \in \Dm}$, $(\X_\dm, \mathcal{X}_\dm)_{\dm \in \Dm} \leftarrow (\X_\dm, \mathcal{X}_\dm)_{\dm \in \Dm}$, $(\w_\dm)_{\dm \in \Dm} \leftarrow (\w_\dm)_{\dm \in \Dm}$, $(f_\dm)_{\dm \in \Dm} \leftarrow (\X_\dm \times \R^A \ni (x,r) \mapsto \delta_\dm \max_{a \in A} \left\{ g_\dm(x,a) + r(a) \right\} \in \R)_{\dm \in \Dm}$, $(\mathcal{F}_\dm^\theta)_{\dm \in \Dm, \; \theta \in \Theta} \leftarrow (\mathcal{F}_\dm^\theta)_{\dm \in \Dm, \; \theta \in \Theta}$, $(X_\dm^\theta)_{\dm \in \Dm, \; \theta \in \Theta} \leftarrow (X_\dm^\theta)_{\dm \in \Dm, \; \theta \in \Theta}$, $(\mathfrak{C}_{n, \dm})_{n \in \N_0, \; \dm \in \Dm} \leftarrow (\mathfrak{C}_{n, \dm})_{n \in \N_0, \; \dm \in \Dm}$, $(V_{n, \dm}^\theta)_{n \in \N_0, \; \dm \in \Dm, \; \theta \in \Theta} \leftarrow (R_{n, \dm}^\theta)_{n \in \N_0, \; \dm \in \Dm, \; \theta \in \Theta}$ in the notation of \cref{mlp_approx_main_theo}) yield the that
    \begin{enumerate}[(1)]
        \item for every $\dm \in \Dm$ there exists a unique function $R_\dm \colon \X_\dm \times A \rightarrow \R$ which is $(\mathcal{X}_\dm \otimes 2^A) / \mathcal{B}(\R)$-measurable \hspace{0.05cm} and \hspace{0.05cm} satisfies \hspace{0.05cm} for \hspace{0.05cm} all \hspace{0.2cm} $(x,a) \in \X_\dm \times A$ \hspace{0.2cm} that \hspace{0.28cm} $\sup_{y \in \X_\dm} \frac{\max_{ b \in A}  | R_\dm(y, b) | }{\w_\dm (y)} < \infty$, $\E \big[ \left| \max_{b \in A} \left\{ g_\dm (X_\dm^{0, x, a}, b) + R_\dm\big( X_\dm^{0,x,a},b \big) \right\} \right| \big] < \infty$,
        \begin{align}\label{eq:mlp_approx_cor_ctrl_R_fp}
            R_\dm(x,a) = \delta_\dm\E \big[ \max_{b \in A} \left\{ g_\dm (X_\dm^{0, x, a}, b) + R_\dm\big( X_\dm^{0,x,a},b \big) \right\}  \big],
        \end{align}
        and
        \item there exist $N \colon (0, 1] \rightarrow \N$ and $c \in \R$ such that for all $\dm \in \Dm$, $\varepsilon \in (0,1]$ it holds that $\mathfrak{C}_{N_\varepsilon , \dm} \leq c \mathfrak{R}_\dm \varepsilon^{-c}$ and
        \begin{align} \label{eq:mlp_approx_cor_ctrl_R_error} 
            \sup_{x \in \X_\dm} \bigg( \frac{\E \big[ \max_{a \in A} | R_\dm(x,a) - R^0_{N_\varepsilon, \dm}(x,a)|^2 \big]}{|\w_\dm(x)|^2} \bigg)^\frac{1}{2} \leq \varepsilon.
        \end{align}
    \end{enumerate}
    For every $\dm \in \Dm$ let $Q_\dm \colon \X_\dm \times A \rightarrow \R$ satisfy for all $(x,a) \in \X_\dm \times A$ that $Q_\dm (x,a) = g_\dm(x,a) + R_\dm(x,a)$. This and item $(1)$ ensure for all $\dm \in \Dm$, $(x,a)\in \X_\dm \times A$ that
    \begin{align}
        Q_\dm(x,a) &= g_\dm(x,a) + R_\dm(x,a) = g_\dm(x,a) + \delta_\dm \E \left[ \max_{b\in A} \left\{ g_\dm(X_\dm^{0,x,a}, b) + R_\dm(X_\dm^{0,x,a}, b)  \right\} \right] \nonumber \\
        &= g_\dm(x, a) + \delta_\dm \E \left[ \max_{b \in A} \left\{ Q_\dm(X_\dm^{0,x,a}, b)\right\} \right]. 
    \end{align}
    Moreover, \hspace{0.1cm} note \hspace{0.1cm} that \hspace{0.1cm} the assumption \hspace{0.1cm} that \hspace{0.1cm} for \hspace{0.1cm} all \hspace{0.1cm} $\dm \in \Dm$, \hspace{0.1cm} $x \in \X_\dm$ \hspace{0.1cm} it \hspace{0.1cm} holds \hspace{0.1cm} that $\max_{a \in A} \left| g_\dm (x,a) \right| \le \kappa \w_\dm(x)$ and item $(1)$ demonstrate that for all $\dm \in \Dm$, $(x,a) \in \X_\dm \times A$ it holds that $\sup_{y \in \X_\dm} \frac{\max_{b \in A} \left| Q_\dm(y,b) \right|}{\w_\dm(y)} < \infty$ and $\E \left[\left| \max_{b \in A} \left\{ Q_\dm(X_\dm^{0,x,a}, b) \right\} \right| \right] < \infty$. Furthermore, for every $\dm \in \Dm$ let $S_\dm \colon \X_\dm \times A \rightarrow \R$ be $(\mathcal{X}_\dm \otimes 2^A) / \mathcal{B}(\R)$-measurable and satisfy for all $\dm \in \Dm$, $(x,a) \in \X_\dm \times A$ that $\sup_{y \in \X_\dm} \frac{\max_{b \in A} \left| S_\dm(y,b) \right|}{\w_\dm(y)} < \infty$, $\E \left[\left| \max_{b \in A} \left\{ S_\dm(X_\dm^{0,x,a}, b) \right\} \right| \right] < \infty$, and $S_\dm(x,a) = g_\dm(x,a) + \delta_\dm \E \left[ \max_{b \in A} \left\{ S_\dm(X_\dm^{0,x,a}, b) \right\} \right]$. It holds for all $\dm \in \Dm$, $(x,a) \in \X_\dm \times A$ that
    \begin{align}
        S_\dm(x,a) - g_\dm(x,a) &= \delta_\dm \E \left[ \max_{b \in A} \left\{ S_\dm(X_\dm^{0,x,a}, b) \right\} \right] \nonumber \\ 
        &= \delta_\dm \E \left[ \max_{b \in A} \left\{ g_\dm(X_\dm^{0,x,a}, b) + S_\dm(X_\dm^{0,x,a}, b) - g_\dm(X_\dm^{0,x,a}, b) \right\} \right].
    \end{align}
    Item $(1)$ implies that for all $\dm \in \Dm$ it holds that $S_\dm - g_\dm = R_\dm$. This establishes for all $\dm \in \Dm$ that $S_\dm = Q_\dm$. This proves item $(i)$. Note that induction and \eqref{eq:mlp_approx_cor_ctrl_R_scheme} demonstrate that for all $\dm \in \Dm$, $n \in \N_0$, $(x,a) \in \X_\dm \times A$, $\theta \in \Theta$ it holds that $Q_{n, \dm}^\theta(x, a) = g_\dm(x, a) + R_{n, \dm}^\theta(x,a)$. Combining this and item $(2)$ implies for all $\dm \in \Dm$, $\varepsilon \in (0,1]$ that
    \begin{align}
        \sup_{x \in \X_\dm} \bigg( \tfrac{\E \big[ \max_{a \in A} | Q_\dm(x,a) - Q_{N_\varepsilon, \dm}^0(x, a) |^2 \big]}{|\w_\dm(x)|^2} \bigg)^\frac{1}{2} = \sup_{x \in \X_\dm} \bigg( \tfrac{\E \big[ \max_{a \in A}| R_\dm(x,a) - R_{N_\varepsilon, \dm}^0(x, a) |^2 \big]}{|\w_\dm(x)|^2} \bigg)^\frac{1}{2} \leq \varepsilon.
    \end{align}
    This proves item $(ii)$. The proof of \cref{mlp_approx_cor_ctrl} is thus completed.
\end{mproof}

\subsection{MLFP approximations for Bellman equations of optimal stopping problems}

\begin{corollary} \label{mlp_approx_cor_stop_slim}
    Let $M \in \N$, let $\Theta = \bigcup_{n \in \N} \Z^n$, let $(\Omega, \mathcal{F}, \mathbb{P})$ be a probability space, let $\Dm$ be a nonempty set, let $\delta_\dm, \mathfrak{R}_\dm \in [0, \infty)$, $\dm \in \Dm$, let $(\X_\dm, \mathcal{X}_\dm)$, $\dm \in \Dm$, be nonempty Borel spaces, for every $\dm \in \Dm$ let $g_\dm \colon \X_\dm \rightarrow \R$ and $G_\dm \colon \X_\dm \rightarrow \R$ be $\mathcal{X}_\dm / \mathcal{B}(\R)$-measurable, for every $\dm \in \Dm$ let $(\mathcal{F}_\dm^\theta)_{\theta \in \Theta}$ be independent sub-$\sigma$-algebras of  $\mathcal{F}$, for every $\dm \in \Dm$ let $X_\dm^\theta = (X_\dm^{\theta, x}(\omega))_{x \in \X_\dm,\; \omega \in \Omega} \colon \X_\dm \times \Omega \rightarrow \X_\dm$, $\theta \in \Theta$, be i.i.d.\ random fields which satisfy for all $\dm \in \Dm$, $\theta \in \Theta$ that $X_\dm^\theta$ is $(\mathcal{X}_\dm \otimes \mathcal{F}_\dm^\theta) / \mathcal{X}_\dm$-measurable, assume $\sup_{\dm \in \Dm} \delta_\dm < 1$, assume $M > \frac{(1 + 3 (\sup_{\dm \in \Dm} \delta_\dm))^2}{(1 - (\sup_{\dm \in \Dm} \delta_\dm))^2}$, assume $\sup_{\dm \in \Dm} \left( \sup_{x \in \X_\dm} |g_\dm(x)| + |G_\dm(x)| \right) < \infty$, for every $\dm \in \Dm$ let $Q_{n, \dm}^\theta \colon \X_\dm \times \Omega \rightarrow \R$, $n \in \N_0$, $\theta \in \Theta$, satisfy for all $n \in \N_0$, $x \in \X_\dm$, $\theta \in \Theta$ that
    \begin{align} \label{mlp_approx_cor_stop_slim_scheme} 
        Q_{n, \dm}^\theta (x) &= g_\dm(x) + \sum_{l = 0}^{n-1} \frac{\delta_\dm}{M^{n-l}} \sum_{i = 1}^{M^{n-l}} \max \left\{ G_\dm \big( X_\dm^{(\theta, l, i), x} \big), Q_{l, \dm}^{(\theta, l, i)}\big(X_\dm^{(\theta, l, i), x} \big) \right\}\\
        &\hspace{4.5cm}- \mathbbm{1}_{\N}(l) \max \left\{ G_\dm \big( X_\dm^{(\theta, l, i), x} \big), Q_{\max \{ l-1, 0 \}, \dm }^{(\theta, -l, i)} \big( X_\dm^{(\theta, l, i),x} \big) \right\} \nonumber
    \end{align}
    and let $\mathfrak{C}_{n, \dm} \in [0, \infty)$, $n \in \N_0$, $\dm \in \Dm$, satisfy for all $n \in \N_0$, $\dm \in \Dm$ that
    \begin{align}
        \mathfrak{C}_{n, \dm} \leq \sum_{l = 0}^{n-1} M^{n-l} \left( \mathfrak{R}_\dm + \mathfrak{C}_{l, \dm} + \mathbbm{1}_\N(l) \mathfrak{C}_{\max \{ l-1, 0 \}, \dm} \right).
    \end{align}
    Then the following holds:
    \begin{enumerate}[$(i)$]
        \item For every $\dm \in \Dm$ there exists a unique function $Q_\dm \colon \X_\dm \rightarrow \R$ which is $\mathcal{X}_\dm/\mathcal{B}(\R)$-measurable and satisfies for all $x \in \X_\dm$ that $\sup_{y \in \X_\dm} |Q_\dm(y)| < \infty$, and
        \begin{align}\label{eq:Bellman_stopping}
            Q_\dm(x) = g_\dm(x) + \delta_\dm \E \left[ \max \left\{ G_\dm(X_\dm^{0,x}), Q_\dm (X_\dm^{0,x}) \right\} \right].
        \end{align}
        \item There exist $N \colon (0,1] \rightarrow \N$ and $c\in \R$ such that for all $\dm \in \Dm$, $\varepsilon \in (0,1]$ it holds that $\mathfrak{C}_{N_\varepsilon, \dm} \leq c \mathfrak{R}_\dm \varepsilon^{-c}$ and
        \begin{align}
            \sup_{x \in \X_\dm} \left( \E \left[ \big| Q_\dm(x) - Q_{N_\varepsilon, \dm}^0(x) \big|^2 \right] \right)^\frac{1}{2} \leq \varepsilon.
        \end{align}
    \end{enumerate}
\end{corollary}

\begin{mproof}{\cref{mlp_approx_cor_stop_slim}}
    Let $\Y_\dm$, $\dm \in \Dm$, be nonempty sets which satisfy that there exists $\Upsilon \in \bigcap_{\dm \in \Dm} \Y_\dm$, and that for all $\dm \in \Dm$ it holds that $\Y_\dm \setminus \{ \Upsilon \} = \X_\dm$, for every $\dm \in \Dm$ let $\mathcal{Y}_\dm = \sigma_{\Y_\dm}(\mathcal{X}_\dm)$, for every $\dm \in \Dm$, $\theta \in \Theta$ let $Y^\theta_\dm = (Y_\dm^{\theta, y, a}(\omega))_{y \in \Y_\dm, \; a \in \{ 0,1 \},\; \omega \in \Omega} \colon \Y_\dm \times \{0,1\} \times \Omega \rightarrow \Y_\dm$ satisfy for all $\dm \in \Dm$, $\theta \in \Theta$, $(y, a) \in \Y_\dm \times \{0,1\}$, $\omega \in \Omega$ that
    \begin{align} \label{mlp_approx_cor_stop_slim_ext_rf}
        Y_\dm^{\theta, y, a}(\omega) &=
        \begin{cases}
            \Upsilon & \colon (y = \Upsilon) \vee (a = 0),\\
            X_\dm^{\theta, y}(\omega) & \colon (y \neq \Upsilon) \wedge (a = 1).
        \end{cases}
    \end{align}
    Note that for all $\dm \in \Dm$ it holds that $Y_\dm^\theta$, $\theta \in \Theta$, are i.i.d.\ random fields and that for every $\dm \in \Dm$, $\theta \in \Theta$ it holds that $Y_\dm^\theta$ is $( \mathcal{Y}_\dm \otimes 2^{\{0,1\}} \otimes \mathcal{F}_\dm^\theta ) / \mathcal{Y}_\dm$-measurable. For every $\dm \in \Dm$ let $h_\dm \colon \Y_\dm \times \{ 0,1\} \rightarrow \R$ satisfy for all $\dm \in \Dm$, $(y, a) \in \Y_\dm \times \{ 0,1 \}$ that
    \begin{align} \label{mlp_approx_cor_stop_slim_ext_payoff}
        h_\dm(y, a) &= 
        \begin{cases}
            0 & \colon y = \Upsilon,\\
            g_\dm(y) & \colon y \neq \Upsilon, \; a = 1,\\
            G_\dm(y) & \colon y \neq \Upsilon, \; a = 0.
        \end{cases}
    \end{align}
    Note that for all $\dm \in \Dm$ it holds that $h_\dm$ is $(\mathcal{Y}_\dm \otimes 2^{\{0,1 \}})/\mathcal{B}(\R)$-measurable. Moreover, for all $\dm \in \Dm$, $(y, a) \in \X_\dm \times \{ 0,1 \}$ it holds that $| h_\dm (y,a) | \leq \max\{  |g_\dm(y)|, |G_\dm(y)| \} \leq \sup_{u \in \Dm} \left( \sup_{x \in \X_\dm} |g_u(x)| + |G_u(x)| \right)$. For every $\dm \in \Dm$ let $q_{n, \dm}^\theta \colon \Y_\dm \times \{ 0,1\} \times \Omega \rightarrow \Y_\dm$, $\theta \in \Theta$, $n \in \N_0$, satisfy for all $n \in \N_0$, $\theta \in \Theta$, $(y, a) \in \Y_\dm \times \{ 0,1\}$ that
    \begin{align} \label{mlp_approx_cor_stop_slim_ext_scheme}
        q_{n, \dm}^\theta(y, a) = h_\dm(y, a) + \sum_{l = 0}^{n-1} \frac{\delta_\dm}{M^{n-l}} \sum_{i = 1}^{M^{n-l}} &\max_{b \in \{ 0,1 \} }\big\{ q_{l, \dm}^{(\theta, l, i)}(Y_\dm^{(\theta, l, i), y,a}, b) \big\} \\
        &-\mathbbm{1}_{\N}(l) \max_{b \in \{ 0,1 \} } \big\{ q_{\max\{ l-1, 0 \}, \dm}^{(\theta, -l, i)}(Y_\dm^{(\theta, l, i), y,a},b) \big\}. \nonumber
    \end{align}
    \cref{mlp_approx_cor_ctrl} (applied with $A \leftarrow \{ 0,1\}$, $\kappa \leftarrow \sup_{\dm \in \Dm} \left( \sup_{x \in \X_\dm} |g_\dm(x)| + |G_\dm(x)| \right)$, $g_\dm \leftarrow h_\dm$, $(\X_\dm, \mathcal{X}_\dm) \leftarrow (\Y_\dm, \mathcal{Y}_\dm)$, $Q_{n, \dm}^\theta \leftarrow q_{n,\dm}^\theta$ for $\dm \in \Dm$ in the notation of \cref{mlp_approx_cor_ctrl}) establishes the following:
    \begin{enumerate}[(1)]
        \item For every $\dm \in \Dm$ there exists a unique function $q_\dm \colon \Y_\dm \times \{ 0,1 \} \rightarrow \R$ which is $(\mathcal{Y}_\dm \otimes 2^{\{ 0,1\}}) / \mathcal{B}(\R)$-measurable and satisfies for all $(y, a) \in \Y_\dm \times \{ 0, 1\}$ that $\sup_{(z, b) \in \Y_\dm \times \{ 0,1 \}} |q_\dm(z,b)| < \infty$
        and
        \begin{align} \label{mlp_approx_cor_stop_slim_ext_fp_eq}
            q_\dm(y, a) = h_\dm(y,a) + \delta_\dm \E \big[ \max_{b \in \{ 0,1 \}} q_d(Y_\dm^{0, y, a}, b) \big].
        \end{align}
        \item[(2)] There exist $N \colon (0, 1] \rightarrow \N$ and $c \in \R$ such that for all $\dm \in \Dm$, $\varepsilon \in (0,1]$ it holds that $\mathfrak{C}_{N_\varepsilon , \dm} \leq c \mathfrak{R}_\dm \varepsilon^c$ and
        \begin{align} \label{mlp_approx_cor_stop_slim_ext_error} 
            \sup_{y \in \Y_\dm} \left( \E \big[ \max_{a \in \{ 0,1\} } | q_\dm(y,a) - q_{N_\varepsilon, \dm}^0(y, a) |^2 \big] \right)^\frac{1}{2} \leq \varepsilon.
        \end{align}
    \end{enumerate}
    Note that for all $\dm \in \Dm$, $(y, a) \in \Y_\dm\setminus \{ \Upsilon\} \times \{ 0 , 1\}$ it holds that $q_\dm(\Upsilon, a) = 0$ and $q_\dm(y, 0) = G_\dm(y)$. Moreover, for all $\dm \in \Dm$, $y \in \Y_\dm \setminus \{ \Upsilon \}$ it holds that 
    \begin{align}
        q_\dm(y, 1) &= h_\dm(y,1) + \delta_\dm \E \big[ \max \big\{ q_\dm(Y_\dm^{0,y,1}, 0), q_\dm(Y_\dm^{0,y,1}, 1) \big\} \big]  \nonumber \\
        &= g_\dm(y) + \delta_\dm \E \left[ \max \left\{ G_\dm (X_\dm^{0, y}), q_\dm(X_\dm^{0,y}, 1) \right\} \right].
    \end{align}
    This and item (1) demonstrate for all $\dm \in \Dm$ that there exists a unique function $Q_\dm \colon \X_\dm \rightarrow \R$ which \hspace{0.05cm} is \hspace{0.05cm} $\mathcal{X}_\dm/\mathcal{B}(\R)$-measurable \hspace{0.05cm} and \hspace{0.05cm} satisfies \hspace{0.05cm} for \hspace{0.05cm} all $x \in \X_\dm$ that $\sup_{y \in \X_\dm} |Q_d(y)| < \infty$ and
    \begin{align}
        Q_\dm(x) = g_\dm(x) + \delta_\dm \E \left[ \max \left\{ G_\dm(X_\dm^{0,x}), Q_\dm(X_\dm^{0,x}) \right\} \right].
    \end{align}
    This proves item $(i)$. Furthermore, note that for all $\dm \in \Dm$, $\theta \in \Theta$, $(y, a) \in \Y_\dm\setminus\{\Upsilon\} \times \{0,1\}$ it holds $\mathbb{P}$-a.s.\ that $q_{0, \dm}^\theta (\Upsilon, a) = 0$, $q_{0, \dm}^\theta(y,0) = G_\dm(y)$, and $q_{0,\dm}^\theta(y, 1) = g_\dm(y)$. This and induction yield that for all $\dm \in \Dm$, $\theta \in \Theta$, $n \in \N_0$, $(y, a) \in \Y_\dm \setminus \{ \Upsilon \} \times \{0,1\}$ it holds $\mathbb{P}$-a.s.\ that $q_{n, \dm}^\theta(\Upsilon, a) = 0$ and $q_{n, \dm}^\theta(y,0) = G_\dm(y)$. Combining this, (\ref{mlp_approx_cor_stop_slim_ext_rf}), (\ref{mlp_approx_cor_stop_slim_ext_payoff}), and (\ref{mlp_approx_cor_stop_slim_ext_scheme}) demonstrates that for all $\dm \in \Dm$, $\theta \in \Theta$, $n \in \N_0$, $y \in \Y_\dm \setminus \{ \Upsilon \} = \X_\dm$ it holds $\mathbb{P}$-a.s.\ that
    \begin{align} \label{mlp_approx_cor_stop_slim_ext_scheme_1}
        q_{n, \dm}^\theta (y, 1) &= h_\dm(y,1) + \sum_{l = 0}^{n-1} \frac{\delta_\dm }{M^{n-l}} \sum_{i = 1}^{M^{n-l}} \max_{b \in \{ 0,1 \} } \left\{ q_{l, \dm}^{(\theta, l, i)}(Y_\dm^{(\theta, l, i), y, 1}, b) \right\} \nonumber \\
        &\hspace{5.5cm}-\mathbbm{1}_{\N}(l) \max_{b \in \{ 0,1 \} } \left\{ q_{\max\{ l-1, 0 \}, \dm}^{(\theta, -l, i)}(Y_\dm^{(\theta, l, i), y, 1}, b) \right\} \nonumber \\
        &= g_\dm(y) + \sum_{l = 0}^{n-1} \frac{\delta_\dm }{M^{n-l}} \sum_{i = 1}^{M^{n-l}} \max \left\{ G_\dm(X_\dm^{(\theta, l, i), y}), q_{l, \dm}^{(\theta, l, i)}(X_\dm^{(\theta, l, i),y}, 1) \right\} \\
        &\hspace{5cm}- \mathbbm{1}_{\N}(l) \max \left\{ G_\dm(X_\dm^{(\theta, l, i),y}), q_{ \max \{ l-1, 0 \} , \dm}^{(\theta, -l, i)}(X_\dm^{(\theta, l, i),y}, 1) \right\}. \nonumber
    \end{align}
    This and (\ref{mlp_approx_cor_stop_slim_scheme}) yield for all $\dm \in \Dm$, $\theta \in \Theta$, $n \in \N_0$, $y \in \Y_\dm \setminus \{ \Upsilon \} = \X_\dm$ it holds $\mathbb{P}$-a.s.\ that
    \begin{align}
        q_{n, \dm}^\theta (y, 1) = Q_{n,\dm}^\theta(y).
    \end{align}
    Combining this and (\ref{mlp_approx_cor_stop_slim_ext_error}) implies that for all $\dm \in \Dm$, $\varepsilon \in (0,1]$ it holds that 
    \begin{align}
        \sup_{x \in \X_\dm} \left( \E \Big[ \big| Q_\dm(x) - Q_{N_\varepsilon, \dm}^0(x) \big|^2 \Big] \right)^\frac{1}{2} \leq \sup_{y \in \Y_\dm} \left( \E \Big[ \max_{a \in \{ 0,1\}} \big| q_\dm(y,a) - q_{N_\varepsilon, \dm}^0(y,a) \big|^2 \Big] \right)^\frac{1}{2} \leq \varepsilon.
    \end{align}
    This establishes item $(ii)$. The proof of \cref{mlp_approx_cor_stop_slim} is thus completed.
\end{mproof}

\subsubsection*{Acknowledgements}
This work has been partially funded by the Deutsche Forschungsgemeinschaft (DFG, \mbox{German} Research Foundation) through the research grant KR 5294/2-1. We also gratefully acknowledge the Cluster of Excellence EXC 2044-390685587, Mathematics Münster: Dynamics-Geometry-Structure funded by the Deutsche Forschungsgemeinschaft (DFG, German Research \mbox{Foundation}).

{
\bibliographystyle{acm}
\bibliography{bibfile}
}

%\newpage
%\printbibliography
\end{document}